\documentclass[11pt]{amsart}
\usepackage{amsmath,amssymb,epsfig,color, amsthm, mathrsfs}
\usepackage{wrapfig,lipsum,floatrow,sidecap,comment}
\usepackage{graphicx}
\usepackage{multirow}
\usepackage{hhline}
\usepackage{url}
\usepackage{mathtools}
\usepackage{longtable}
\graphicspath{ {images/} }
\allowdisplaybreaks
\graphicspath{ {images/} }

\def\P{\mathcal{P}}
\def\p{\prime}

\newtheorem{theorem}{Theorem}[section]

\newtheorem{corollary}[theorem]{Corollary}
\newtheorem{definition}[theorem]{Definition}

\newtheorem{remark}[theorem]{Remark}
\newtheorem{example}[theorem]{Example}
\begin{document}
\title[Enumeration of Pretzel Links]{The Enumeration of Alternating Oriented Pretzel Links}
\author{Charlotte Aspinwall$^\flat$, Tobias Clark$^\sharp$ and Yuanan Diao$^\dagger$}
\address{$^\flat$ Department of Mathematics\\
Cornell University\\
Ithaca, NY 14850}
\email{csa74@cornell.edu}
\address{$^\sharp$ Department of Mathematics\\
Vassar College\\
Poughkeepsie, NY 12604}
\email{tclark@vassar.edu}
\address{$^\dagger$ Department of Mathematics and Statistics\\
University of North Carolina Charlotte\\
Charlotte, NC 28223}
\email{ydiao@charlotte.edu}
\thanks{$^\flat, ^\sharp$: supported by NSF-REU DMS-2150179.}

\begin{abstract}
In this paper, we tabulate the set of alternating oriented pretzel links. We derive a closed formula for the precise number of alternating oriented pretzel links with any given crossing number $c$. Numerical computation suggests that this number grows approximately at the rate of $0.155e^{0.588c}$.    
\end{abstract} 

\maketitle

\section{Introduction}

Ever since the beginning of knot theory, knot tabulation, a subject of study first started by Tait in the 1870s, has been a focus of study.
This is a very difficult problem. Tait (with Little and Kirkman) managed to create the first knot table containing prime knots with up to 10 crossings entirely by hand, and it only contained one duplicate (which was discovered by Perko in 1973 \cite{Perko1974}). Fast forwarding to the late 20th century, armed with more modern techniques and computing powers, Hoste, Thistlethwaite and Weeks \cite{Hoste1998} tabulated all prime knots with up to 16 crossings (there are $1,701,936$ of them). Recently, Burton \cite{Burton2020} extended the tabulation of all prime knots with up to 19 crossings (and there are 350 millions of them) and Thistlethwaite \cite{Thistlethwaite} tabulated all prime knots with 20 crossings (there are more than 1.8 billions of them). The tabulation of links with multiple components is even more difficult and much less is known.

\medskip
Instead of trying to tabulate all knots/links with a fixed crossing number, one can tabulate specific classes of knots and links with a given crossing number. For some knot/link classes such as the torus knots/links, there are relatively few of them for each given crossing number, hence the tabulation is easier. The class of all rational knots and links is the first {\em large class} of prime knots/links (oriented or non-oriented) that has been completely tabulated \cite{Diao2022, Ernst1987}. That is, for any given crossing number $c$, the numbers of all oriented/non-oriented rational knots and rational links have been precisely determined. 
By the term``large class" we mean that the number of rational knots and links with crossing number $c$ grows exponentially as a function of $c$. 
In this paper, we seek to extend the table to include another large class of (oriented) links, namely the set of all (oriented) alternating oriented pretzel links. 
We should note that rational links and pretzel links are both subclasses of Montesinos links. 
While the rational links possess some very nice and unique properties (which make them the most studied large class of links), the pretzel links provide a good supplement as they possess certain properties not shared by rational links. 
For example, all rational links are alternating but most pretzel links are not. 
Another example is that all rational links have bridge index 2. 
This means that they are simple links as measured by the bridge index while pretzel links can have arbitrarily large bridge indices \cite{Diao2003}. 
As our study shows, the class of pretzel links is also a large one, so it certainly makes sense for us to have a better understanding of this class. 

\medskip
In our approach, we divide the alternating oriented pretzel links into three different Types: Type 1, Type 2 and Type 3 alternating oriented pretzel links. 
This is because the counting for each Type requires a different tool. 
For each given crossing number $c$ and Type of alternating oriented pretzel link, we are able to obtain the precise number of links in the form of a closed formula. The formulations for Type 1 and Type 2 are relatively simple, but the formulation for Type 3 alternating oriented pretzel links is quite a bit more complicated. 
This is due to the more complicated structure of these pretzel links. 
Numerical results based on these formulations suggest that the number of oriented alternating oriented pretzel links with crossing number $c$, grows approximately as $0.155e^{0.588c}$, indicating that the number of alternating oriented pretzel links with a given crossing number $c$ grows exponentially in terms of $c$.  

\medskip
We will organize the rest of this paper as follows. In the next section, we shall introduce the basic concepts, notations and terminology used throughout the paper. The three Types of alternating oriented pretzel links will be defined. In Section \ref{sec_class}, we establish the foundation for our tabulation, namely the if and only if conditions that allow us to determine whether two oriented alternating oriented pretzel link diagrams represent the same link. In Section \ref{results}, we present our main results. In Section \ref{sec_enum}, we enumerate the oriented alternating oriented pretzel links for each of the three Types and prove the theorems presented in Section \ref{results}. Finally, in Section \ref{sec_num}, we present the numerical results based on the theoretical formulas outlined in Section \ref{sec_notation}.

\section{Notations, terminology and preliminary results}\label{sec_notation}
\begin{figure}[htb!]
\includegraphics[scale=0.9]{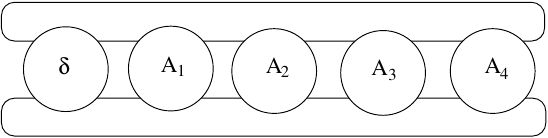}
\caption{A diagram depicting a general Montesinos link with $4$ rational tangles and $\delta$ horizontal half-twists (crossings).}
\label{Montesinos}
\end{figure}

A Montesinos link is a link that can be represented by a diagram as shown in Figure \ref{Montesinos} where the first circle contains a horizontal row of $|\delta|\ge 0$ half twists (crossings) while each other circle $A_i$ for $1 \leq i \leq k$ represents a two string rational tangle with at least 2 crossings. Figure \ref{Montesinos} depicts the case $k=4$. Each rational tangle $A_i$ can be represented by a rational number $p_i/q_i$ such that $p_i$ and $q_i$ are co-prime and $q_i>1$. Similar to the notation used in \cite{Burde}, we can denote a Montesinos link diagram by $\mathbb{M}(\delta;p1/q_1,\ldots,p_k/q_k)$. Notice that in general a Montesinos link may not be alternating and the notation $\mathbb{M}(\delta;p1/q_1,\ldots,p_k/q_k)$ is used for the un-oriented Montesinos links. 

\medskip
In the case that $k\le 2$, a Montesinos link is just a rational link. 
A pretzel link is a special Montesinos link in which each rational tangle consists of a vertical strip of $2$ or more half twists (crossings). In our case, we are interested in the enumeration of alternating oriented pretzel links with a given crossing number. 
To avoid double counting a rational link as a pretzel link, we restrict ourselves to the pretzel links whose Montesinos diagrams consist of $k\ge 3$ vertical strips each containing $2$ or more crossings. 

\begin{figure}[!hbt]
\begin{center}
\includegraphics[scale=0.8]{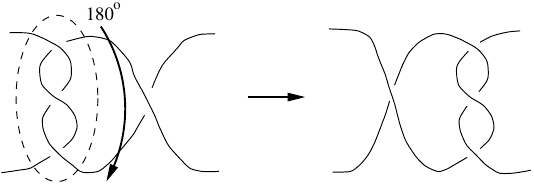}
\end{center}
\caption{A flype move applied to a single crossing and an adjacent strip in a pretzel link diagram.
\label{flypes}}
\end{figure}

\begin{figure}[!hbt]
\begin{center}
\includegraphics[scale=0.7]{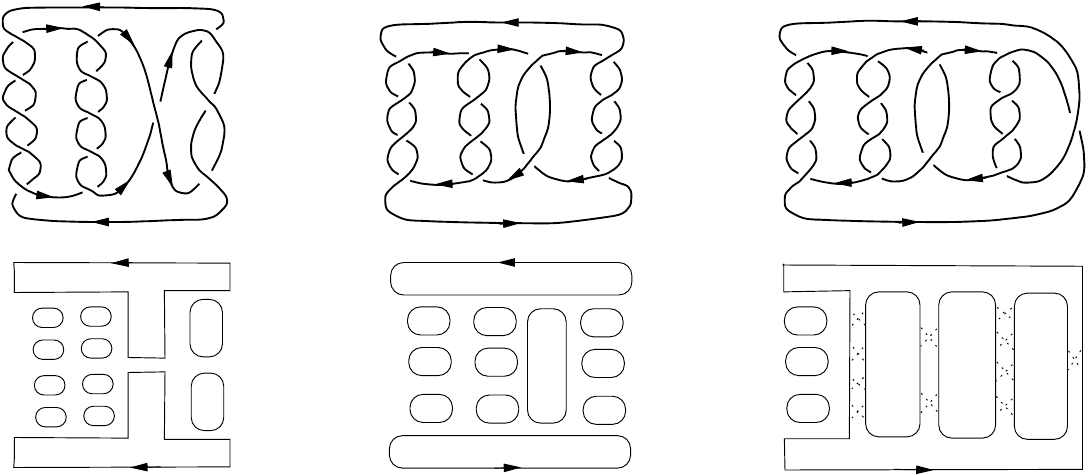}
\end{center}
\caption{Examples of alternating oriented pretzel links.
\label{pretzel_fig}}
\end{figure}

\medskip
A flype, as shown in Figure \ref{flypes}, is a link type preserving move that also preserves the components and their orientations in the case of an oriented link. Using repeated flypes if needed, one can place a crossing in $\delta$ of an alternating oriented link diagram between any two vertical strips of the diagram, not necessarily at the left side of the diagram. That is, diagrams as shown in Figure \ref{pretzel_fig} are all alternating oriented pretzel link diagrams. The following theorem is well known.

\medskip
\begin{theorem}\label{flype_thm} \cite{Menasco1990} Two reduced alternating oriented link diagrams $D_1$ and $D_2$ are {\em equivalent} if and only if there exists a sequence of flypes that take $D_1$ to $D_2$. Here, $D_1$ to $D_2$ are said to be equivalent if there exists an ambient isotopy that takes $D_1$ to $D_2$.
\end{theorem}

\medskip
Each alternating oriented pretzel link diagram naturally corresponds to a sequence of non-zero integers (called a {\em code} of the diagram) such that: 
(i) the absolute value of an integer in the sequence indicates the number of crossings in its corresponding strip, and 
(ii) the sign of the integer indicates the crossing sign of the crossings in its corresponding strip. For example, the three pretzel links in Figure \ref{pretzel_fig} can be coded by $(5,5,1,3)$, $(-4,-4,-2,-4)$ and $(-4,4,2,4,1)$. 

\begin{remark}\label{flype_remark}{\em
Let $D$ be an alternating oriented pretzel link diagram.
We leave it to our reader to verify that a flype that can be performed on $D$ either involves only crossings within a vertical strip of $D$ (or within $\delta$), or involves a crossing in the $\delta$ group and some vertical strips of $D$ (plus possibly some crossings in $\delta$ as well). In the first case, the flype does not change $D$ at all, and in the second case the flype only changes the relative position of the crossing in $\delta$ with the strips involved. The flype does not change the signs of the crossings involved nor how the crossings are smoothed. If $\delta=0$, then there does not exist any flype that would change $D$. If $\delta>0$ and we only insist to place the crossings in $\delta$ adjacent to each other, then these crossings can be placed between any two adjacent vertical strips of $D$. But no other forms of $D$ exist under this restriction.
}
\end{remark}

\medskip
We now divide the set of all alternating oriented pretzel link diagrams into three disjoint subsets, which we shall call Type 1, Type 2 and Type 3 alternating oriented pretzel link diagrams as defined in \cite{Diao2024_1,Diao2024_2}. 
The reason for this division is that different counting techniques are needed in order to count the alternating oriented pretzel link diagrams in these subsets. 

\medskip
\begin{definition}{\em
Let $D$ be an alternating oriented pretzel link diagram and let $S(D)$ be the Seifert circle decomposition of $D$. 
Let $C_1$ and $C_2$ be the Seifert circles containing the top and bottom strands of $D$, respectively. 
We say that $D$ is a Type 1 (2) alternating oriented pretzel link diagram if $C_1\not= C_2$ and $C_1$, $C_2$ have different (same) orientations. 
In the case that $C_1=C_2$, we say that $D$ is a Type 3 alternating oriented pretzel link diagram.
}
\end{definition}

\medskip
For example, the pretzel link diagrams coded by $(5,5,1,3)$, $(-4,-4,-2,-4)$ and $(-4,4,2,4,1)$ in the top of Figure \ref{pretzel_fig} are of Types 1, 2 and 3 respectively according to their Seifert circle decompositions as shown at the bottom of Figure \ref{pretzel_fig}. 

\begin{figure}[!hbt]
\begin{center}
\includegraphics[scale=0.7]{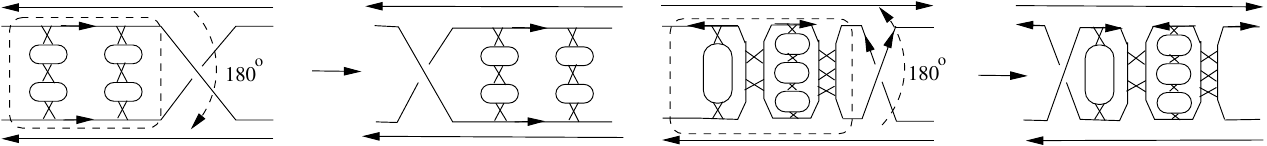}
\end{center}
\caption{The effect of a flype involving a crossing in $\delta$ on the Seifert circle decomposition of a Type 1 (left) or a Type 3 (right) alternating oriented pretzel link diagram.
\label{Type13_flype}}
\end{figure}

\begin{remark}\label{flype_remark2}{\em By Remark \ref{flype_remark}, a Type 2 alternating oriented pretzel link diagram admits no flypes that would change the diagram, and flypes on a Type 1 or 3  alternating oriented pretzel link diagram, under the restriction that all crossings in $\delta$ are placed adjacent to each other, can only change the crossings in $\delta$ are placed, not the signs of the crossings and how they are smoothed. That is, Type 1 (Type 3) alternating oriented pretzel link diagrams will remain Type 1 (Type 3) alternating oriented pretzel link diagrams after flypes. (This is easy to see from the perspective of the Seifert circle decomposition of the diagram as shown in Figure \ref{Type13_flype}.) Consequently, alternating oriented pretzel link diagrams with different Types are never equivalent by Theorem \ref{flype_thm}. }
\end{remark}

\begin{remark}\label{Type_smooth_remark}{\em
Let $D$ be an alternating oriented pretzel link diagram. We point out a few more facts about $D$ below, but leave the verifications of them to the reader.\\
\noindent
(1) If the crossings in $\delta$ smooth horizontally (vertically), then $D$ is of  Type 1 (Type 3); \\
(2) If $D$ is of Type 1, then each strip contains an odd number of crossings;\\
\noindent
(3) If $D$ is of Type 2, then each strip contains an even number of crossings;;\\
\noindent
(4) If $D$ is of Type 1 or Type 2, then the crossings in $D$ are either all positive or all negative;\\
\noindent
(5) If $D$ is of Type 3, then the positive crossings will be smoothed either all horizontally or all vertically, while the negative crossings are smoothed the opposite way;\\
\noindent
(6) If $D$ is of Type 3, then the number of strips whose crossings are to be smoothed vertically must be even if the horizontal row of $\delta$ crossings are counted as $\delta$ vertical strips (since these crossings will also be smoothed vertically);\\
\noindent
(7) If $D$ is of Type 3 and the crossings in a strip of $D$ are to be smoothed horizontally, then the number of crossings in this strip must be even;\\
\noindent
}
\end{remark}

\medskip
\begin{remark}\label{mirror_remark}{\em
By Theorem \ref{flype_thm}, Remark \ref{flype_remark2}, and (4), (5) of Remark \ref{Type_smooth_remark}, the mirror image of an alternating oriented pretzel link diagram $D$ is never equivalent to $D$ itself. This is relatively easy to see if $D$ is of Type 1 or 2, since all crossings of $D$ would be positive or negative, while the crossings of its mirror image have opposite sign. If $D$ is of Type 3 and the crossings in $D$ that are smoothed vertically are all positive (negative), then the mirror images of these crossings are precisely the ones that are smoothed in the mirror image of $D$, but they are all negative (positive). The conclusion now follows from Remark \ref{flype_remark2}. }
\end{remark}

\medskip
Since our purpose is to count the number of distinct alternating oriented pretzel links with a given crossing number, and each alternating oriented pretzel link and its mirror image are always distinct, we only need to count one of them, and simply double our total count at the end. For this reason, in the rest of the paper, we shall only count Type 1 and Type 2 alternating oriented pretzel links with positive crossings, and Type 3 alternating oriented pretzel links in which positive crossings are smoothed vertically and negative crossings are smoothed horizontally.

\medskip
Without loss of generality, we can always assume that the bottom long strand of the diagram is oriented from right to left (since otherwise we can rotate the diagram around the $y$-axis by $180^\circ$ to achieve this). Throughout the rest of the paper, we shall use the term {\em standard diagram} to denote an alternating oriented pretzel link diagram with the bottom long strand of the diagram oriented from right to left and with all crossings in $\delta$ placed at the left side of the diagram. By Remark \ref{Type_smooth_remark}, a Type 1 standard diagram with positive crossings and $k\ge 3$ tangles must have a code of the form $(\delta;2\alpha_1+1,2\alpha_2+1,\ldots, 2\alpha_k+1)$ with $\delta\ge 0$ and $\alpha_j\ge 1, 1\le j\le k$. 
However this code by itself does not tell us the diagram is of Type 1 or of Type 3. For example one can construct a Type 1 standard diagram and a Type 3 standard diagram both with the same code $(2;3,3,3,3)$: in the first construction every crossing is smoothed horizontally and in the second construction every crossing is smoothed vertically. Thus, to avoid this ambiguity we use the notation 
$P_1(\delta;2\alpha_1+1,2\alpha_2+1,\ldots, 2\alpha_k+1)$ to indicate a Type 1 standard diagram with positive crossings. 
Similarly, we shall use
$P_2(2\alpha_1,2\alpha_2,\ldots, 2\alpha_k)$ ($\alpha_j\ge 1, 1\le j\le k$) to denote a standard Type 2 alternating link diagram with positive crossings. 
Finally, we shall use the notation
$P_3(\delta; \beta_1,\beta_2, \ldots, \beta_k)$ to indicate a Type 3 standard diagram in which positive crossings are smoothed vertically and negative crossings are smoothed horizontally.  In this case, it is implied that $k\ge 3$ and $|\beta_j |\ge 2, 1\le j\le k$. 
We shall call these $P_1$, $P_2$ and $P_3$ notations the {\em $t$-codes} of the standard diagrams.

\medskip
\begin{example}{\em
The $t$-codes corresponding to the standard diagrams of the alternating oriented pretzel links in Figure \ref{pretzel_fig} are $P_1(1;5,5,3)$, $P_2(-4,-4,-2,-4)$ and $P_3(1;-4,4,2,4)$ respectively.
}
\end{example}

\begin{remark}\label{standard_diag_remark}{\em 
There is a unique standard diagram corresponding to each $t$-code. 
As a consequence of this, from now on, we shall treat a $t$-code the same as the unique standard diagram it corresponds to. 
That is, we shall say that two $t$-codes are equivalent if their corresponding standard diagrams represent equivalent oriented links. Notice that every alternating oriented pretzel link diagram has at least one standard diagram, hence corresponds to at least one $t$-code. However, two different $t$-codes may correspond to two alternating oriented pretzel links that are equivalent, in which case we shall say that the two $t$-codes are {\em equivalent}.
Throughout the paper, we shall use the notation $\cong$ to denote the equivalency of two alternating oriented pretzel links, as well as the equivalency of two $t$-codes.
}
\end{remark}

The statement in the first sentence of Remark \ref{standard_diag_remark} requires a proof.
Let us consider a Type 1 standard diagram $D$ with positive crossings. 
The other two types can be similarly discussed and are left to the reader. 
The top long strand of $D$ is also oriented from right to left since $D$ is of Type 1, so the top strand must be the over-strand at the top crossing in the left most strip of the diagram since the crossing is positive and is smoothed horizontally. 
Now there is only one way to draw the other crossings in the diagram since the diagram is alternating. 
 
\section{Classification of the standard diagrams}\label{sec_class}

By Remark \ref{standard_diag_remark}, every alternating oriented pretzel link can be represented by at least one $t$-code. Thus, if we can determine the equivalency of $t$-codes, we would be able to count the number of distinct alternating oriented pretzel links by counting the number of equivalent classes of $t$-codes. A classification theorem for the $t$-codes is thus in order. By Remarks \ref{flype_remark2} and \ref{mirror_remark}, we only need to consider the equivalency of Type 1 and Type 2 $t$-codes with positive entries, and Type 3 $t$-codes in which the crossings in $\delta$ are positive (if there are any) and the positive crossings are the ones that are smoothed vertically. 

\medskip
\begin{theorem}\label{T2}
The following statements hold:\\
\noindent
(i) $P_1(\delta;2\alpha_1+1,2\alpha_2+1,\ldots, 2\alpha_k+1)\cong P_1(\delta^\prime;2\alpha^\prime_1+1,2\alpha^\prime_2+1,\ldots, 2\alpha^\prime_{\ell}+1)$ if and only $\delta=\delta^\prime$, $k=\ell$, and $(2\alpha^\prime_1+1,\ldots, 2\alpha^\prime_{k}+1)$ is a cyclic permutation of $(2\alpha_1+1,2\alpha_2+1,\ldots, 2\alpha_k+1)$.
\\
\noindent
(ii) $P_2(2\alpha_1,2\alpha_2,\ldots, 2\alpha_k) \cong P_2(2\alpha^\prime_1,2\alpha^\prime_2,\ldots, 2\alpha^\prime_{\ell})$ if and only if $k=\ell$ and $(2\alpha^\prime_1,2\alpha^\prime_2,\ldots, 2\alpha^\prime_{k})$ is either a cyclic permutation of $(2\alpha_1,2\alpha_2,\ldots, 2\alpha_k)$, or a cyclic permutation of $(2\alpha_k,\ldots,2\alpha_2, 2\alpha_1)$.
\\
\noindent
(iii) $P_3(\delta; \beta_1,\beta_2, \ldots, \beta_k)\cong P_3(\delta^\prime; \beta^\prime_1,\beta^\prime_2, \ldots, \beta^\prime_{\ell})$ if and only if $\delta=\delta^\prime$, $k=\ell$, and $(\beta^\prime_1,\beta^\prime_2, \ldots, \beta^\prime_{k})$ is either a cyclic permutation of $( \beta_1,\beta_2, \ldots, \beta_k)$, or a cyclic permutation of $(\beta_k,\ldots,\beta_2,  \beta_1)$.
\end{theorem}

\begin{proof}
By Theorem \ref{flype_thm} and Remark \ref{flype_remark2}, if $D_1$ and $D_2$ are two standard alternating oriented pretzel link diagrams that are equivalent, then $D_2$ can only be obtained from $D_1$ through a combination of the following moves: (a) flypes that change the positions of the crossings in $\delta$; (b) a cyclic permutation of the vertical strips in $D_1$; (c) a rotation of the diagram around the $x$-axis by $180^\circ$ and (d) a rotation of the diagram around the $y$-axis by $180^\circ$. Among these moves, the flypes are the only moves that actually alter the diagram, the other moves merely provide views of the same diagram from different perspectives. 

\medskip
(i) In the case of a Type 1 standard diagram, a rotation around the $x$-axis by $180^\circ$ does not change the diagram and a rotation of the diagram around the $y$-axis by $180^\circ$ does not result in a standard diagram since the bottom strand is not oriented from right to left. Figure \ref{Type1_eqv} demonstrates that a cyclic permutation of the vertical strips of a Type 1 standard diagram combined with suitable flypes leads to an equivalent Type 1 standard diagram:  the first three figures show the standard diagram (with $k=4$) from a cyclic view point. The fourth figure shows the cyclic view of the standard diagram with $t$-code $P_1(\delta; 2\alpha_2+1,\ldots, 2\alpha_k+1, 2\alpha_1+1)$, which is obtained by flypes that move the crossings in $\delta$ to between the strips $2\alpha_1+1$ and $2\alpha_2+1$, and so on. In other words, any $t$-code equivalent to $P_1(\delta;2\alpha_1+1,2\alpha_2+1,\ldots, 2\alpha_k+1)$ is the $t$-code of a Type 1 standard diagram obtained from the diagram of $P_1(\delta;2\alpha_1+1,2\alpha_2+1,\ldots, 2\alpha_k+1)$ by a cyclic permutation of $(2\alpha_1+1,2\alpha_2+1,\ldots, 2\alpha_k+1)$ followed by flypes that move the crossings in $\delta$ to the left side of the diagram.

\begin{figure}[!hbt]
\begin{center}
\includegraphics[scale=0.7]{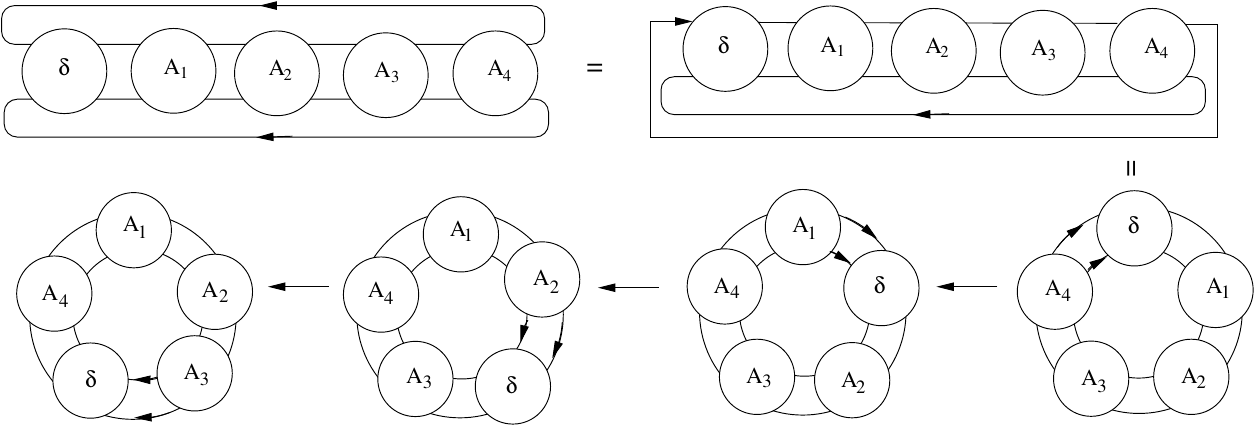}
\end{center}
\caption{A cyclic view point of the Type 1 standard diagrams after using flypes to place the crossings in $\delta$ between two strips. Shown is the case of $k=4$ with $A_j$ being a vertical strip containing $2\alpha_j+1$ crossings. \label{Type1_eqv}}
\end{figure}

\smallskip
(ii) In the case of a Type 2 standard diagram, a rotation around the $x$-axis by $180^\circ$ alone does not yield a standard diagram, but this rotation followed by another rotation  around the $y$-axis by $180^\circ$ results in a Type 2 standard diagram. Cyclic permutations of the strips in these two diagrams provide the other different perspectives of the same diagram since there are no flypes that can alter the diagram. That is, if $P_2(2\alpha_1,2\alpha_2,\ldots, 2\alpha_k) \cong P_2(2\alpha^\prime_1,2\alpha^\prime_2,\ldots, 2\alpha^\prime_{\ell})$, then $(2\alpha^\prime_1,2\alpha^\prime_2,\ldots, 2\alpha^\prime_{\ell})$ must be either a cyclic permutation of $(2\alpha_1,2\alpha_2,\ldots, 2\alpha_k)$, or a cyclic permutation of $(2\alpha_k,2\alpha_{k-1},\ldots, 2\alpha_1)$ (which implies that $\ell=k$).

\smallskip
(iii) In the case of a Type 3 standard diagram, a rotation around the $x$-axis by $180^\circ$ alone does not yield a standard diagram, but this rotation followed by another rotation  around the $y$-axis by $180^\circ$ will result in a Type 3 standard diagram. Thus all the other equivalent standard diagrams can be obtained from these two by cyclic permutations of the vertical strips combined with suitable flypes keeping the crossings in $\delta$ at the left side of the diagram. 
That is, if $P_3(\delta; \beta_1,\beta_2, \ldots,\beta_k)\cong P_3(\delta^\prime; \beta^\prime_1,\beta^\prime_2, \ldots, \beta^\prime_{\ell})$, then $\delta=\delta^\prime$ and $(\beta^\prime_1,\beta^\prime_2,\ldots,\beta^\prime_{\ell})$ is either a cyclic permutation of $(\beta_1,\beta_2,\ldots,\beta_k)$, or a cyclic permutation of $(\beta_k,\beta_{k-1},\ldots,\beta_1)$ (which of course also implies that $\ell=k$). \end{proof}

\section{The main results}\label{results}

 \begin{theorem}\label{Type1_Thm}
Let $\P_1(c)$ be the number of Type 1 alternating oriented pretzel links  with crossing number $c$ and with positive crossings.  Then we have
\begin{equation}\label{P1equation}
\P_1(c)=\sum N(\frac{c-\delta-k}{2},k),
\end{equation}
where the summation is taken over all $\delta\ge 0$, $k\ge 3$, $n\ge 3k+\delta$ such that $n-\delta-k$ is even, and 
\begin{equation}\label{N(n,k)equation}
N(n,k)=\frac{1}{k}\sum_{d | gcd(n,k)} \phi(d)\binom{\frac{n}{d}-1}{\frac{k}{d}-1}
\end{equation}
with $\phi(d)$ being the Euler totient function.
\end{theorem}

\begin{theorem}\label{Type2_Thm}
Let $\P_2(c)$ be the number of Type 2 alternating oriented pretzel links with crossing number $c$  and with positive crossings.  Then we have
\begin{equation}\label{Type2_eq}
\P_2(c)=
\begin{cases}
\sum_{3\le k\le n}B(n,k),&\ {\rm if}\ c=2n\ge 6;\\
0,&\ {\rm otherwise}
\end{cases}
\end{equation}
where 
\begin{equation}\label{B(n,k)eq_1}
B(n,k)=\frac{1}{2}N(n,k)+h(n,k)
\end{equation}
with
\begin{equation}\label{h(n,k)eq}
h(n,k)=
\begin{cases}
\frac{1}{2}\binom{\frac{n-1}{2}}{\frac{k-1}{2}} & \ {\rm if}\ n\equiv k\equiv 1\ {\rm mod}\ 2,\\
\frac{1}{2}\binom{\frac{n}{2}-1}{\frac{k-1}{2}} & \ {\rm if}\ n\equiv 0\ {\rm mod}\ 2,\ k\equiv 1\ {\rm mod}\ 2,\\
\frac{1}{2}\binom{\frac{n-1}{2}}{\frac{k}{2}} & \ {\rm if}\ n\equiv 1\ {\rm mod}\ 2,\ k\equiv 0\ {\rm mod}\ 2,\\
\frac{1}{2}\binom{\frac{n}{2}}{\frac{k}{2}} & \ {\rm if}\ n\equiv 0\ {\rm mod}\ 2,\ k\equiv 0\ {\rm mod}\ 2.\\
\end{cases}
\end{equation}
\end{theorem}

\begin{theorem}\label{Type3_Thm}
Consider the set of all Type 3 alternating oriented pretzel links whose positive crossings are the ones that are smoothed vertically.
Let $\P_3(c)$ be the number of such Type 3 alternating oriented pretzel links with crossing number $c$.  Then we have
\begin{equation}\label{P3(c)equ}
\P_3(c)=\sum_{(\delta;n_1,k_1;n_2,k_2)\in \Omega_c}B(n_1,k_1;n_2,k_2),
\end{equation}
where 
\begin{equation}\label{Type3_bracelet_formula}
B(n_1,k_1;n_2,k_2)=
\frac{1}{2k}\sum_{d| gcd(k_1,k_2,n_1,n_2)}\phi(d)\binom{k/d}{k_1/d}\binom{n_1/d-1}{k_1/d-1}\binom{n_2/d-1}{k_2/d-1}+h(n_1,k_1;n_2,k_2),
\end{equation}
 $k=k_1+k_2$, $h(n_1,k_1;n_2,k_2)$ is as defined in formulas (\ref{hformula1}) to (\ref{hformula3}), and $\Omega_c$ is the set of quintuples of non-negative integers $(\delta;n_1,k_1;n_2,k_2)$ satisfying the following conditions:  (i) $n_1\ge k_1\ge 0$ and $n_2\ge k_2\ge 0$; (ii) $\delta+k_1+n_1+2n_2=c$; (iii) $\delta+k_1\ge 2$ and must be even; (iv) $k_1+k_2\ge 3$.
\end{theorem}

\medskip
Formula (\ref{hformula1}) applies when $k=k_1+k_2$ is odd:
\begin{equation}\label{hformula1}
h(n_1,k_1;n_2,k_2)
=
\begin{cases}
 \frac{1}{2}\binom{\frac{k-1}{2}}{\frac{k_2}{2}}\binom{\frac{n_1-1}{2}}{\frac{k_1-1}{2}}\binom{\frac{n_2}{2}-1}{\frac{k_2}{2}-1}, & \ {\rm if}\ n_1\equiv k_1\equiv 1\ {\rm mod}\ 2,\ n_2\equiv k_2\equiv 0 \ {\rm mod}\ 2,\\
 \frac{1}{2}\binom{\frac{k-1}{2}}{\frac{k_1}{2}}\binom{\frac{n_1}{2}-1}{\frac{k_1}{2}-1}\binom{\frac{n_2-1}{2}}{\frac{k_2-1}{2}}, & \ {\rm if}\ n_1\equiv k_1\equiv 0\ {\rm mod}\ 2,\ n_2\equiv k_2\equiv 1 \ {\rm mod}\ 2,\\
 \frac{1}{2}\binom{\frac{k-1}{2}}{\frac{k_1}{2}}\binom{\frac{n_1}{2}-1}{\frac{k_1}{2}-1}\binom{\frac{n_2}{2}-1}{\frac{k_2-1}{2}}, & \ {\rm if}\ n_1\equiv k_1\equiv n_2\equiv 0\ {\rm mod}\ 2,\ k_2\equiv 1 \ {\rm mod}\ 2,\\
  \frac{1}{2}\binom{\frac{k-1}{2}}{\frac{k_2}{2}}\binom{\frac{n_1}{2}-1}{\frac{k_1-1}{2}}\binom{\frac{n_2}{2}-1}{\frac{k_2}{2}-1}, & \ {\rm if}\ n_1\equiv k_2\equiv n_2\equiv 0\ {\rm mod}\ 2,\ k_1\equiv 1 \ {\rm mod}\ 2,\\
 0, & \ {\rm if}\  k_1\equiv  n_2\equiv  1\ {\rm mod}\ 2\ {\rm or}\  k_2\equiv  n_1\equiv  1\ {\rm mod}\ 2.
  \end{cases}
  \end{equation}
Formula (\ref{hformula2}) applies when $k$ is even but $k_1$, $k_2$ are both odd:
\begin{equation}\label{hformula2}
h(n_1,k_1;n_2,k_2)
=
\begin{cases}
 \frac{1}{2}\binom{\frac{k}{2}-1}{\frac{k_1-1}{2}}\binom{\frac{n_1-1}{2}}{\frac{k_1-1}{2}}\binom{\frac{n_2-1}{2}}{\frac{k_2-1}{2}}, & \ {\rm if}\ n_1\equiv n_2\equiv 1\ {\rm mod}\ 2,\\
 \frac{1}{2}\binom{\frac{k}{2}-1}{\frac{k_1-1}{2}}\binom{\frac{n_1}{2}-1}{\frac{k_1-1}{2}}\binom{\frac{n_2-1}{2}}{\frac{k_2-1}{2}}, & \ {\rm if}\ n_1\equiv n_2+1\equiv 0\ {\rm mod}\ 2,\\
 \frac{1}{2}\binom{\frac{k}{2}-1}{\frac{k_1-1}{2}}\binom{\frac{n_1-1}{2}}{\frac{k_1-1}{2}}\binom{\frac{n_2}{2}-1}{\frac{k_2-1}{2}}, & \ {\rm if}\ n_1\equiv n_2+1\equiv 1\ {\rm mod}\ 2,\\
 \frac{1}{2}\binom{\frac{k}{2}-1}{\frac{k_1-1}{2}}\binom{\frac{n_1}{2}-1}{\frac{k_1-1}{2}}\binom{\frac{n_2}{2}-1}{\frac{k_2-1}{2}}, & \ {\rm if}\ n_1\equiv n_2\equiv 0\ {\rm mod}\ 2.
\end{cases}
\end{equation}
Formula (\ref{hformula3}) applies when $k$, $k_1$ and $k_2$ are all even:
\begin{equation}\label{hformula3}
h(n_1,k_1;n_2,k_2)
=
\begin{cases}
 \phantom{+}\frac{1}{4}\binom{\frac{k}{2}-1}{\frac{k_2}{2}}\binom{\frac{n_2}{2}-1}{\frac{k_2}{2}-1}\left(\binom{\frac{n_1}{2}}{\frac{k_1}{2}}+\binom{\frac{n_1}{2}-1}{\frac{k_1}{2}}\right) & \\
 +  \frac{1}{4}\binom{\frac{k}{2}-1}{\frac{k_1}{2}}\binom{\frac{n_1}{2}-1}{\frac{k_1}{2}-1}
 \left(\binom{\frac{n_2}{2}}{\frac{k_2}{2}}+\binom{\frac{n_2}{2}-1}{\frac{k_2}{2}}\right) &  \\
 + \frac{1}{4}\binom{\frac{k}{2}}{\frac{k_1}{2}}\binom{\frac{n_1}{2}-1}{\frac{k_1}{2}-1}\binom{\frac{n_2}{2}-1}{\frac{k_2}{2}-1}, & \ {\rm if}\ n_1\equiv n_2\equiv 0\ {\rm mod}\ 2,\\
\frac{1}{2}\binom{\frac{k}{2}-1}{\frac{k_1}{2}}\binom{\frac{n_1}{2}-1}{\frac{k_1}{2}-1}\binom{\frac{n_2-1}{2}}{\frac{k_2}{2}}, & \ {\rm if}\ n_1\equiv n_2+1\equiv 0\ {\rm mod}\ 2,\\ 
\frac{1}{2}\binom{\frac{k}{2}-1}{\frac{k_2}{2}}\binom{\frac{n_2}{2}-1}{\frac{k_2}{2}-1}\binom{\frac{n_1-1}{2}}{\frac{k_1}{2}}, & \ {\rm if}\ n_1+1\equiv n_2\equiv 0\ {\rm mod}\ 2,\\
0, &\ {\rm if}\ n_1\equiv n_2\equiv 1\ {\rm mod}\ 2.\\
\end{cases}
\end{equation}

\medskip
Let us remind our reader that by Remark \ref{mirror_remark}, the total number of alternating oriented pretzel links with crossing number $c$ is thus given by $2(\P_1(c)+\P_2(c)+\P_3(c))$.

\section{The proofs of Theorems \ref{Type1_Thm}, \ref{Type2_Thm} and \ref{Type3_Thm}}\label{sec_enum}

By Remark \ref{standard_diag_remark}, the enumeration of alternating oriented pretzel links with a given crossing number $c$ is just an enumeration of the equivalent classes of $t$-codes whose (absolute) parameters sum to $c$ with the equivalence relations given by by Theorem \ref{T2} (i)--(iii). 
In the next three subsections, tools and known results based on cycle index functions and generating functions in classic combinatorics will be used. 
For the details, we shall refer our reader to a textbook covering these topics such as \cite{West2021}. 

\subsection{The derivation of Type 1 enumeration formula}

\medskip
Let us consider the set of all $t$-codes $P_1(\delta; 2\alpha_1+1,\ldots,2\alpha_k+1)$ with a fixed $\delta\ge 0$ and $k\ge 3$ (keep in mind that $\alpha_j\ge 1$ for each $j$) under the condition $\delta+k+2\sum_{j=1}^k\alpha_j=c$. 
By Theorem \ref{T2}(i), two such $t$-codes $P_1(\delta; 2\alpha_1+1,\ldots,2\alpha_k+1)$ and $P_1(\delta; 2\alpha^\p_1+1,\ldots,2\alpha^\p_k+1)$ are equivalent if and only if $(2\alpha^\p_1+1,\ldots,2\alpha^\p_k+1)$ is a cyclic permutation of $(2\alpha_1+1,\ldots,2\alpha_k+1)$, or if and only if $(\alpha^\p_1,\ldots,\alpha^\p_k)$ is a cyclic permutation of $(\alpha_1,\ldots,\alpha_k)$. Counting the number of equivalence classes of the $k$-tuples of positive integer entries where the equivalence is defined by cyclic permutation is a classical combinatorics problem. 
In our case, the $k$-tuples are subject to an additional condition $\sum_{1\le j\le k}\alpha_j=n$ for some integer $n\ge k$, since our goal here is to enumerate Type 1 alternating oriented pretzel links with a fixed crossing number. 
So, let $X(n,k)$ be the set of all $k$-tuples of positive integer entries subject to the additional condition that the sum of the entries is equal to $n$. 
Let $N(n,k)$ be the number of cyclic equivalence classes of the $k$-tuples in $X(n,k)$. 
$N(n,k)$ can be interpreted as the number of necklaces with $k$ beads, where each bead is colored by a positive integer. 
This same integer is the ``weight" of the color and it must be the case that the total weight of the weights of the colors is $n$. 
Under this setting, the generating function for the color set is simply $f(t)=t+t^2+t^3+\cdots=\frac{t}{1-t}$, and the P\'{o}lya-Redfield Theorem (\cite{Polya1937,Polya1987,Redfield1927}) asserts that the generating function for $N(n,k)$ can be obtained from the cycle index function of the cyclic group $C_k$:
$$
Z_{C_k}(x_1,\ldots,x_k)=\frac{1}{k}\sum_{d| k}\phi(d)x_d^{\frac{k}{d}}
$$
(where $\phi(d)$ is the Euler totient function) by replacing $x_d$ with $f(t^d)$. 
Thus the generating function for $N(n,k)$ is
$$
N_k(t)=\frac{1}{k}\sum_{d | k} \phi(d)\left(\frac{t^d}{1-t^d}\right)^{\frac{k}{d}}=\frac{1}{k}\sum_{d | k} \phi(d)t^k\left(\frac{1}{1-t^d}\right)^{\frac{k}{d}}
$$
and $N(n,k)$ is the coefficient of $t^n$ in the above summation. 
Substitution $s=t^d$ in the identity
$$
\left(\frac{1}{1-s}\right)^m=\sum_{j\ge 0}\binom{m+j-1}{m-1}s^j,
$$
we then have:
\begin{eqnarray*}
\left(\frac{t^d}{1-t^d}\right)^{\frac{k}{d}}&=&t^k\left(\frac{1}{1-t^d}\right)^{\frac{k}{d}}
=t^k\sum_{j\ge 0}\binom{\frac{k}{d}+j-1}{\frac{k}{d}-1}t^{dj}.
\end{eqnarray*}
If this power series is to contain a $t^n$ term, then we have $k+jd=n$ for some $j$. 
Since $k$ is a multiple of $d$, it means we must have $d|n$ as well, namely $d| gcd(n,k)$, and $j=(n-k)/d$. That is, $\left(\frac{t^d}{1-t^d}\right)^{\frac{k}{d}}$ contains a $t^n$ term if and only if $d| gcd(n,k)$. 
Furthermore, in the case that $d| gcd(n,k)$, the $t^n$ term in $\left(\frac{t^d}{1-t^d}\right)^{\frac{k}{d}}$ is $\binom{\frac{n}{d}-1}{\frac{k}{d}-1}t^{n}$. 
It follows that
\begin{equation}\label{N(n,k)equation}
N(n,k)=\frac{1}{k}\sum_{d | gcd(n,k)} \phi(d)\binom{\frac{n}{d}-1}{\frac{k}{d}-1}.
\end{equation}

\medskip
We can now express $\P_1(c)$ in terms of the function $N(n,k)$. 
For a fixed $\delta\ge 0$ and $k\ge 3$, if $P_1(\delta; 2\alpha_1+1,\ldots,2\alpha_k+1)$ has crossing number $c$, then we must have
 $$
 c=\delta+k+2\sum_{1\le j\le k}\alpha_j,
 $$ 
 hence $c-\delta-k$ must be even (and must be at least since $\alpha_j\ge 1$ and $k\ge 3$) and
 $$
 \sum_{1\le j\le k}\alpha_j=\frac{c-\delta-k}{2}.
 $$
 We see that $N(\frac{c-\delta-k}{2},k)$ then gives us the precise number of Type 1 alternating oriented pretzel links with crossing number $c$ and parameters $\delta$ and $k$. 
 This leads to 
 $$
\P_1(c)=\sum N(\frac{c-\delta-k}{2},k),
$$
where the summation is taken over all $\delta\ge 0$, $k\ge 3$ satisfying the conditions that $n-\delta-k$ is even and $n-\delta\ge 3k$. 
This proves Theorem \ref{Type1_Thm}.

\medskip
We can simplify the summation in (\ref{P1equation}) by defining $N(n,k)$ to be zero if $n<k$. 
\begin{corollary}\label{Cor1}
If $c=2q+1(\ge 9)$ is odd, then 
\begin{equation}\label{L1e1}
\P_1(c)=\sum_{i\ge 3}\sum_{j\ge \lfloor \frac{i}{2}\rfloor}N(q-j,i),
\end{equation}
On the other hand, if $c=2q(\ge 10)$ is even, then
\begin{equation}\label{L1e2}
\P_1(c)=\sum_{i\ge 3}\sum_{j\ge \lceil \frac{i}{2}\rceil}N(q-j,i),
\end{equation}
\end{corollary}

\begin{proof} Consider the case when $c$ is odd.
For any given even $i$, $\delta$ must be odd, hence $\delta+i=2j+1$ with $j=\frac{\delta-1}{2}+\frac{i}{2}\ge \frac{i}{2}=\lfloor \frac{i}{2}\rfloor$. Similarly, if $i$ is odd, then $\delta\ge 0$ must be even. 
Again $\delta+i=2j+1$ must be odd and $j=\frac{\delta}{2}+\frac{i-1}{2}\ge \frac{i-1}{2}=\lfloor \frac{i}{2}\rfloor$. 
The formula now follows. 
The case of $c$ even can be similarly discussed and is left to the reader.
\end{proof}

\subsection{The derivation of Type 2 enumeration formula} We now consider the set of all $t$-codes $P_2(2\alpha_1,\ldots,2\alpha_k)$ with a fixed $\delta\ge 0$ and $k\ge 3$ (keep in mind that $\alpha_j\ge 1$ for each $j$) under the condition $\delta+k+2\sum_{j=1}^k\alpha_j=c$. 
By Theorem \ref{T2}(ii), two $t$-codes $P_2(2\alpha_1,\ldots,2\alpha_k)$ and $P_2(2\alpha^\p_1,\ldots,2\alpha^\p_k)$ are equivalent if and only if $(2\alpha^\p_1,\ldots,2\alpha^\p_k)$ is either a cyclic permutation of $(2\alpha_1,\ldots,2\alpha_k)$, or a cyclic permutation of $(2\alpha_k,\ldots,2\alpha_1)$. This can be rephrased as: $P_2(2\alpha_1,\ldots,2\alpha_k)\cong P_2(2\alpha^\p_1,\ldots,2\alpha^\p_k)$ if and only if $(\alpha^\p_1,\ldots,\alpha^\p_k)$ can be obtained from $(\alpha_1,\ldots,\alpha_k)$ by performing a dihedral group $D_k$ action.

\medskip
We thus define $B(n,k)$ similarly to $N(n,k)$, except that the equivalence class of a $k$-tuple now includes the reversal of the $k$-tuple and its cyclic permutations. 
These equivalent classes are the orbits of the dihedral group $D_k$ (which is of order $2k$) acting on the set of all $k$-tuples with positive integer entries whose sum equals $n$. 
The cycle index function for $D_k$ is
\begin{equation}\label{cycle_index_dihedral}
Z_{D_{k}}(x_1,\ldots,x_k)=
\begin{cases}
\frac{1}{2k}\sum_{d|k} \phi(d)x_d^{k/d} +\frac{1}{2} x_1x_2^{(k-1)/2}   &\mbox{for odd $k$,}\\
\frac{1}{2k}\sum_{d|k} \phi(d)x_d^{n/d}
+\frac{1}{4} x_1^2x_2^{(k-2)/2}+\frac{1}{4} x_2^{k/2}   &\mbox{for even $k$.}\\
\end{cases}  
\end{equation}
As in the case of $N(n,k)$, by the P\'{o}lya-Redfield Theorem, the generating function for $B(n,k)$ is 
\begin{eqnarray}
B_k(t)&=&
Z_{D_k}(f(t),f(t^2),\cdots, f(t^k))\\
&=&\frac{1}{2}N_k(t)+\begin{cases}
\frac{1}{2} \frac{t}{1-t}\left(\frac{t^2}{1-t^2}\right)^{(k-1)/2}  &\mbox{for odd $k$,}\\
\frac{1}{4}\left(\frac{t}{1-t}\right)^2\left(\frac{t^2}{1-t^2}\right)^{(k-2)/2}+\frac{1}{4} \left(\frac{t^2}{1-t^2}\right)^{k/2}
 &\mbox{for even $k$}.\nonumber\\
\end{cases}\label{B_kt}
\end{eqnarray}
That is, $B(n,k)$ is the coefficient of the $t^n$ power term in $B_k(t)$. Direct computation leads to 
\begin{equation}\label{B(n,k)eq_1}
B(n,k)=\frac{1}{2k}\sum_{d | gcd(n,k)} \phi(d)\binom{\frac{n}{d}-1}{\frac{k}{d}-1}+h(n,k)=\frac{1}{2}N(n,k)+h(n,k)
\end{equation}
where 
\begin{equation}\label{h(n,k)eq}
h(n,k)=
\begin{cases}
\frac{1}{2}\binom{\frac{n-1}{2}}{\frac{k-1}{2}} & \ {\rm if}\ n\equiv k\equiv 1\ {\rm mod}\ 2,\\
\frac{1}{2}\binom{\frac{n}{2}-1}{\frac{k-1}{2}} & \ {\rm if}\ n\equiv 0\ {\rm mod}\ 2,\ k\equiv 1\ {\rm mod}\ 2,\\
\frac{1}{2}\binom{\frac{n-1}{2}}{\frac{k}{2}} & \ {\rm if}\ n\equiv 1\ {\rm mod}\ 2,\ k\equiv 0\ {\rm mod}\ 2,\\
\frac{1}{2}\binom{\frac{n}{2}}{\frac{k}{2}} & \ {\rm if}\ n\equiv 0\ {\rm mod}\ 2,\ k\equiv 0\ {\rm mod}\ 2.\\
\end{cases}
\end{equation}
We note that the formulation for $h(n,k)$ was derived more than a century ago by Sommerville \cite{Sommerville}.

\medskip
Let $c$ be a given crossing number and $\P_2(c)$ the number of mirror image pairs of Type 2 alternating oriented pretzel links with $c$ crossings and at least 3 strips. 
Keep in mind that here each strip contains at least $2$ crossings and must contain an even number of crossings. 
Obviously, if $c$ is odd, or $c<5$, then $\P_2(c)=0$. If $c=2n\ge 6$ is even, then a Type 2 alternating oriented pretzel link has a $t$-code of the form $P_2(2\alpha_1,2\alpha_2,\ldots,2\alpha_k)$ such that $\alpha_j\ge 1$ and $\sum_{1\le j\le k}2\alpha_j=c$, or $\sum_{1\le j\le k}\alpha_j=n$. 
That is, each such Type 2 alternating oriented pretzel link corresponds to a $k$-tuple $(\alpha_1,\alpha_2,\ldots,\alpha_k)$ with $\sum_{1\le j\le k}\alpha_j=n$. 
Therefore, there are precisely $B(n,k)$ Type 2 alternating oriented pretzel links with positive crossings and $k\ge 3$ strips. 
That is,
$$
\P_2(c)=
\begin{cases}
\sum_{3\le k\le n}B(n,k),&\ {\rm if}\ c=2n\ge 6;\\
0&\ {\rm otherwise.}
\end{cases}
$$
This proves Theorem \ref{Type2_Thm}.

\begin{example}{\em If $c=14$, then $n=14/2=7$ and the total number of Type 2 alternating oriented pretzel links with $14$ crossings (that are all positive) is
\begin{eqnarray*}
\P_2(14)&=&\sum_{3\le k\le 7}B(7,k)\\
&=&B(7,3)+B(7,4)+B(7,5)+B(7,6)+B(7,7)\\
&=&\frac{\frac{1}{3}\binom{6}{2}+\binom{3}{1}}{2}+\frac{\frac{1}{4}\binom{6}{3}+\binom{3}{2}}{2}+\frac{\frac{1}{5}\binom{6}{4}+\binom{3}{2}}{2}\\
&+&\frac{\frac{1}{6}\binom{6}{5}+\binom{3}{3}}{2}+\frac{\frac{1}{7}(\binom{6}{6}+6\binom{0}{0})+\binom{3}{3}}{2}\\
&=&
4+4+3+1+1=13.
\end{eqnarray*}
Direct construction also shows that there are 13 such Type 2 alternating oriented pretzel links with positive crossings. 
}
\end{example}

\subsection{The derivation of Type 3 enumeration formula}
Let $c$ be a given crossing number and let $\P_3(c)$ be the number of Type 3 alternating oriented pretzel links with the following properties: (a) their crossing numbers are all $c$; (b) they contain at least 3 vertical strips and (c) the positive crossings are smoothed vertically (hence if $\delta>0$ then the crossings in $\delta$ are also positive). Let $\Gamma(c)$ be the set of all $t$-codes corresponding to standard Type 3 alternating oriented pretzel link diagrams satisfying conditions (a), (b) and (c). So $\P_3(c)$ is just the number of equivalence classes of $t$-codes in $\Gamma(c)$. 
Let us further divide $\Gamma(c)$ into subsets $\Gamma_c(\delta;m_1,k_1;m_2,k_2)$ such that $t$-codes in $\Gamma_c(\delta;m_1,k_1;m_2,k_2)$ share  common parameters $\delta$, $k_1$, $k_2$, $m_1$ and $m_2$, where $k_1$ ($k_2$) and $m_1$ ($m_2$) are the numbers of positive (negative) strips and the total numbers of crossings in the positive (negative) strips respectively. 
For example, the parameters of $P_3(3;3,-4,-2,4,-6,2,5,-6,7)$ are $c=42$, $\delta=3$, $k_1=5$, $m_1=21$, $k_2=4$ and $m_3=18$. That is, $P_3(3;3,-4,-2,4,-6,2,5,-6,7)\in \Gamma_{42}(3;21,5;18,4)$.

\medskip
By Theorem \ref{T2}(iii), $\Gamma_c(\delta;m_1,k_1;m_2,k_2)\cap \Gamma_c(\delta^\p;m^\p_1,k^\p_1;m^\p_2,k^\p_2)=\emptyset$ if $(\delta^\p;m^\p_1,k^\p_1;m^\p_2,k^\p_2)$ is not identical to $(\delta;m_1,k_1;m_2,k_2)$. 
Thus we can obtain $\P_3(c)$ by summing the numbers of equivalent $t$-codes from each non-empty $\Gamma_c(\delta;m_1,k_1;m_2,k_2)$. 
Notice that $\Gamma_c(\delta;m_1,k_1;m_2,k_2)$ is non-empty if and only if the following conditions all hold: $m_1\ge 2k_1$, $m_2\ge 2k_2$, $m_2$ is even, $k_1+k_2\ge 3$, $\delta+k_1\ge 2$ is even and $\delta+m_1+m_2=c$.

\medskip
By subtracting one from each positive $\beta_j$ and dividing each negative $\beta_j$ by 2, we note that each $t$-code $P_3(\delta;\beta_1,\beta_2,\ldots,\beta_k)$ in $\Gamma_c(\delta;m_1,k_1;m_2,k_2)$ corresponds uniquely to a $k$-tuple of non-zero integers $(\alpha_1,\alpha_2,...,\alpha_k)$ with a parameter quadruple $(n_1,k_1;n_2,k_2)$ satisfying the following conditions: 
(i) $n_1(=m_1-k_1)\ge k_1\ge 0$ and $n_2(=m_2/2)\ge k_2\ge 0$; 
(ii) $\delta+k_1+n_1+2n_2(=\delta+m_1+m_2)=c$; 
(iii) $\delta+k_1\ge 2$ and must be even; 
(iv) $k_1+k_2\ge 3$. 
We shall denote by $\Omega_c$ the set of all quintuples $(\delta;n_1,k_1;n_2,k_2)$ satisfying conditions (i)--(iv) above.
By Theorem \ref{T2}(iii), two $t$-codes with the same parameter quintuple $(\delta;m_1,k_1;m_2,k_2)$ represent the same Type 3 alternating oriented pretzel link if and only if their corresponding $k$-tuples with parameter quintuple $(\delta;n_1,k_1;n_2,k_2)$ are related by a Dihedral group $D_k$ action, where $k=k_1+k_2\ge 3$, $n_1=m_1-k_1\ge k_1\ge 0$ and $n_2=m_2/2\ge k_2\ge 0$. 
That is, the problem of counting the number of Type 3 alternating oriented pretzel links represented by $t$-codes in the set $\Gamma_c(\delta;m_1,k_1;m_2,k_2)$ is the same as counting the equivalence classes of $k$-tuples with parameter quintuples $(\delta;n_1,k_1;n_2,k_2)$ in $\Omega_c$ where the equivalence classes are defined by the orbits of Dihedral group $D_k$ actions.

\medskip
Thus we have the same cycle index function 
$Z_{D_{k}}(x_1,\ldots,x_k)$ as defined in (\ref{cycle_index_dihedral}). 
However, a $k$-tuple with parameters $(n_1,k_1;n_2,k_2)\in \Omega_c$ is a bracelet whose beads are colored from two different sets of colors: those colored by a positive integer $\alpha_j$ with weight $\alpha_j$, and those colored by a negative integer $\alpha_j$ with weight $|\alpha_j|$. 
The total weight of the positively colored beads is $n_1$ and the total weight of the negatively colored beads is $n_2$. 
Thus, in this case, we can choose our generating function as 
$$
f(t)+g(s)=t+t^2+\cdots +s+s^2+\cdots=\frac{t}{1-t}+\frac{s}{1-s}
$$
where the powers of the variable $t$ correspond to the positive entries in the $k$-tuples and the powers of the variable $s$ correspond to the negative entries in the $k$-tuples. 
Let $B(n_1,k_1; n_2,k_2)$ be the number of bracelets with $k_1$ ($k_2$) being the number of positive (negative) entries in the $k$-tuple, and $n_1$ ($n_2$) be the total weight of the positive (negative) entries. 
By the P\'{o}lya-Redfield Theorem, the generating function for $B(n_1,k_1; n_2,k_2)$ is
\begin{eqnarray}
B_{k}(t,s)&=&
Z_{D_k}(f(t)+g(s),f(t^2)+g(s^2),\cdots, f(t^k)+g(s^k))\label{B_{k}(t,s)}\label{B_{k}(t,s)}\\
&=&\frac{1}{2k}\sum_{d | k} \phi(d)\left(\frac{t^d}{1-t^d}+\frac{s^d}{1-s^d}\right)^{\frac{k}{d}}\nonumber\\
&+&
\begin{cases}
\frac{1}{2} \left(\frac{t}{1-t}+\frac{s}{1-s}\right)\left(\frac{t^2}{1-t^2}+\frac{s^2}{1-s^2}\right)^{(k-1)/2}   &\mbox{for odd $k$,}\nonumber\\
\frac{1}{4} \left(\frac{t}{1-t}+\frac{s}{1-s}\right)^2\left(\frac{t^2}{1-t^2}+\frac{s^2}{1-s^2}\right)^{(k-2)/2}+\frac{1}{4} \left(\frac{t^2}{1-t^2}+\frac{s^2}{1-s^2}\right)^{k/2}   &\mbox{for even $k$}.\nonumber\\
\end{cases}
\end{eqnarray}
Notice that $B(n_1,k_1; n_2,k_2)$ is the coefficient of the $t^{n_1}s^{n_2}$ power term in $B_{k}(t,s)$, which must be contributed by ``$k_1$ copies of $f(t)$ and $k_2$ copies of $f(s)$". 
Here $f^2(t)$ and $f(t^2)$ would both count as 2 copies of $f(t)$, $f^3(t)$, $f(t)f(t^2)$, and $f(t^3)$ etc would all count as 3 copies of $f(t)$, and so on. 
Direct computations then lead to equation (\ref{Type3_bracelet_formula}), that is,
\begin{equation}\label{B(n_1,k_1;n_2,k_2)}
B(n_1,k_1;n_2,k_2)=
\frac{1}{2k}\sum_{d| gcd(k_1,k_2,n_1,n_2)}\phi(d)\binom{k/d}{k_1/d}\binom{n_1/d-1}{k_1/d-1}\binom{n_2/d-1}{k_2/d-1}+h(n_1,k_1;n_2,k_2)
\end{equation}
where $h(n_1,k_1;n_2,k_2)$ is as given in formulas (\ref{hformula1}) to (\ref{hformula3}).
Of course, in the case that $k_2=0$ or $k_1=0$, we must have $n_2=0$ or $n_1=0$ respectively, in which case $B(n_1,k_1;n_2,k_2)$ simply becomes $B(n_1,k_1)$ (if $k_2=n_2=0$) or $B(n_2,k_2)$ (if $k_1=n_1=0$). This leads to the following formulation of $\P_3(c)$ 
$$
\P_3(c)=\sum_{(\delta;n_1,k_1;n_2,k_2)\in \Omega_c}B(n_1,k_1;n_2,k_2)
$$
as stated in Theorem \ref{Type3_Thm}. 
This concludes the proof of Theorem \ref{Type3_Thm}.

\medskip
In the following we shall demonstrate how our computations are done for the derivation of 
formula (\ref{Type3_bracelet_formula}) using $h(n_1,k_1;n_2,k_2)$ for the case when $n_1$, $k_1$, $n_2$ and $k_2$ are all even. 
The other parts of the formula can be similarly derived and in fact are slightly easier, and are placed in the appendix section.

\medskip
First let us examine the term 
$$
\frac{1}{4} \left(f(t)+f(s)\right)^2\left(f(t^2)+f(s^2)\right)^{\frac{k-2}{2}}=\frac{1}{4} \left(\frac{t}{1-t}+\frac{s}{1-s}\right)^2\left(\frac{t^2}{1-t^2}+\frac{s^2}{1-s^2}\right)^{\frac{k-2}{2}}.
$$
Here we do not have to consider the term $2f(t)f(s)$ in $(f(t)+f(s))^2=f^2(t)+2f(t)f(s)+f^2(s)$ since every term in $f(t)f(s)(f(t^2)+f(s^2))^{\frac{k-2}{2}}$ is contributed by an odd number of copies of $f(t)$ and $f(s)$ but we only need to consider the terms contributed by even number of copies of $f(t)$ and $f(s)$.
That is, we only need to consider 
$$
\frac{1}{4} \left((\frac{t}{1-t})^2+(\frac{s}{1-s})^2\right)\left(\frac{t^2}{1-t^2}+\frac{s^2}{1-s^2}\right)^{(k-2)/2}.
$$
Let us first try to extract the $t^{n_1}s^{n_2}$ term from 
$$
\frac{1}{4} \left(\frac{t}{1-t}\right)^2\left(\frac{t^2}{1-t^2}+\frac{s^2}{1-s^2}\right)^{(k-2)/2}
$$
using $k_1$ copies of the generating function $f(t)$ and $k_2$ copies of the generating function $f(s)$. Since $\left(\frac{t}{1-t}\right)^2$ already represents two copies of $f(t)$, we need $k_1-2$ copies of $f(t)$ from $\left(\frac{t^2}{1-t^2}+\frac{s^2}{1-s^2}\right)^{(k-2)/2}$. The term in the binomial expansion of $\left(\frac{t^2}{1-t^2}+\frac{s^2}{1-s^2}\right)^{(k-2)/2}$ with $k_1-2$ copies of $f(t)$ is 
$$
\binom{\frac{k}{2}-1}{\frac{k_2}{2}}\left(\frac{t^2}{1-t^2}\right)^{\frac{k_1}{2}-1}\left(\frac{s^2}{1-s^2}\right)^{\frac{k_2}{2}}.
$$
That is, we need to extract the $t^{n_1}s^{n_2}$ term from
\begin{equation}\label{he1}
\frac{1}{4}\binom{\frac{k}{2}-1}{\frac{k_2}{2}} \left(\frac{t}{1-t}\right)^2\left(\frac{t^2}{1-t^2}\right)^{\frac{k_1}{2}-1}\left(\frac{s^2}{1-s^2}\right)^{\frac{k_2}{2}}.
\end{equation}
Using the Maclaurin power series formula 
\begin{equation}\label{MacSeries}
(1-z)^{-m}=\sum_{j\ge 0}\binom{m+j-1}{j}z^j, m\ge 1,
\end{equation}
 it is relatively easy to see that $s^{n_2}$ is contributed from $\left(\frac{s^2}{1-s^2}\right)^{\frac{k_2}{2}}=s^{k_2}\left(\frac{1}{1-s^2}\right)^{\frac{k_2}{2}}$ and is of the form $\binom{n_2/2-1}{k_2/2-1}s^{n_2}$. For $k_1=2$, the contribution to $t^{n_1}$ comes from $(t/(1-t))^2$, which is $(n_1-1)t^{n_1}$ since 
\begin{eqnarray*}
\left(\frac{t}{1-t}\right)^2&=&
 t^{2}\left(1+2t+3t^2+\cdots\right).
 \end{eqnarray*}

\medskip
For $k_1>2$, $\left(\frac{t}{1-t}\right)^2\left(\frac{t^2}{1-t^2}\right)^{\frac{k_1}{2}-1}$ can be expressed as
\begin{equation}\label{expans_k1even}
t^{k_1}\left(1+2t+3t^2+\cdots\right)\left(1+\binom{\frac{k_1}{2}-2+1}{\frac{k_1}{2}-2}t^2+\cdots+\binom{\frac{k_1}{2}-2+j}{\frac{k_1}{2}-2}t^{2j}+\cdots\right),
\end{equation}
so the coefficient of $t^{n_1}$ from  it is
$$
\sum_{0\le j\le (n_1-k_1)/2}(2j+1)\binom{\frac{n_1}{2}-2-j}{\frac{k_1}{2}-2}.
$$
Using the binomial coefficient identities
 \begin{equation}\label{binom_iden_1}
 \binom{n}{k}+ \binom{n-1}{k}+ \binom{n-2}{k}+\cdots+ \binom{k}{k}= \binom{n+1}{k+1},
 \end{equation}
 \begin{equation}\label{binom_iden_2}
 \binom{n-1}{k}+2 \binom{n-2}{k}+3 \binom{n-3}{k}+\cdots+(n-k) \binom{k}{k}= \binom{n+1}{k+2},
 \end{equation}
 and 
 \begin{equation}\label{binom_iden_3}
 \binom{n}{k-1}+ \binom{n}{k}= \binom{n+1}{k+1},
 \end{equation}
we have
\begin{eqnarray*}
&&\sum_{0\le j\le (n_1-k_1)/2}(2j+1)\binom{\frac{n_1}{2}-2-j}{\frac{k_1}{2}-2}\\
&=&\sum_{0\le j\le (n_1-k_1)/2}\binom{\frac{n_1}{2}-2-j}{\frac{k_1}{2}-2}+
2\sum_{1\le j\le (n_1-k_1)/2}j\binom{\frac{n_1}{2}-2-j}{\frac{k_1}{2}-2}\\
&=&
\binom{\frac{n_1}{2}-1}{\frac{k_1}{2}-1}+2\binom{\frac{n_1}{2}-1}{\frac{k_1}{2}}=\binom{\frac{n_1}{2}}{\frac{k_1}{2}}+\binom{\frac{n_1}{2}-1}{\frac{k_1}{2}}.
 \end{eqnarray*}
 Thus the contribution to the coefficient of $t^{n_1}s^{n_2}$ from (\ref{he1}) is
 \begin{equation}\label{t^{n_1}s^{n_2}_1}
 \frac{1}{4}\binom{\frac{k}{2}-1}{\frac{k_2}{2}}\binom{\frac{n_2}{2}-1}{\frac{k_2}{2}-1}\left(\binom{\frac{n_1}{2}}{\frac{k_1}{2}}+\binom{\frac{n_1}{2}-1}{\frac{k_1}{2}}\right).
 \end{equation}
 Notice that the above formula also works for the case $k_1=2$ and $n_1>2$. If $n_1=2$, then we must also have $k_1=2$ and we still have $\binom{\frac{n_1}{2}}{\frac{k_1}{2}}+\binom{\frac{n_1}{2}-1}{\frac{k_1}{2}}=n_1-1$ since $\binom{0}{1}=0$.
 Similarly, the contribution to the coefficient of $t^{n_1}s^{n_2}$ from 
 $$
\frac{1}{4} \left(\frac{s}{1-s}\right)^2\left(\frac{t^2}{1-t^2}+\frac{s^2}{1-s^2}\right)^{(k-2)/2}
$$
 is
 \begin{equation}\label{t^{n_1}s^{n_2}_2}
 \frac{1}{4}\binom{\frac{k}{2}-1}{\frac{k_1}{2}}\binom{\frac{n_1}{2}-1}{\frac{k_1}{2}-1}
 \left(\binom{\frac{n_2}{2}}{\frac{k_2}{2}}+\binom{\frac{n_2}{2}-1}{\frac{k_2}{2}}\right).
 \end{equation}
 Finally, let us consider the contribution from 
$$
\frac{1}{4} \left(\frac{t^2}{1-t^2}+\frac{s^2}{1-s^2}\right)^{k/2}
$$
which must be from 
$$
\frac{t^{k_1}s^{k_2}}{4}\binom{\frac{k}{2}}{\frac{k_1}{2}} \left(\frac{1}{1-t^2}\right)^{\frac{k_1}{2}}\left(\frac{1}{1-s^2}\right)^{\frac{k_2}{2}}.
$$
 Thus we obtain a coefficient 
 \begin{equation}\label{t^{n_1}s^{n_2}_3}
 \frac{1}{4}\binom{\frac{k}{2}}{\frac{k_1}{2}}\binom{\frac{n_1}{2}-1}{\frac{k_1}{2}-1}\binom{\frac{n_2}{2}-1}{\frac{k_2}{2}-1}.
 \end{equation}
  (\ref{hformula3}) then follows after we combine (\ref{t^{n_1}s^{n_2}_1}), (\ref{t^{n_1}s^{n_2}_2}) and (\ref{t^{n_1}s^{n_2}_3}).

\medskip
\begin{example}{\em
Let us compute $\P_3(10)$ using two approaches: one by formula (\ref{P3(c)equ}) and one by direction calculation. In order to use (\ref{P3(c)equ}), we list the set of all quintuples $(\delta;n_1,k_1;n_2,k_2)$ of the parameters in $\Omega_{10}$. } 

\begin{center}
\begin{tabular}{ccccc}
$(4;2,2;1,1)$&$(4;0,0;3,3)$&$(3;4,3;0,0)$& $(3;2,1;2,2)$& $(2;4,4;0,0)$\\
 $(2;2,2;2,2)$& $(2;4,2;1,1)$& $(2;2,2;2,1)$& $(2;0,0;4,3)$& $(2;0,0;4,4)$\\
$(1;4,3;1,1)$&$(1;6,3;0,0)$& $(1;2,1;3,3)$& $(1;2,1;3,2)$& $(1;4,1;2,2)$\\
 $(0;6,4;0,0)$&$(0;4,4;1,1)$& $(0;4,2;2,2)$& $(0;2,2;3,2)$& $(0;2,2;3,3)$\\
  $(0;6,2;1,1)$& $(0;4,2;2,1)$& $(0;2,2;3,1)$ &&
  \end{tabular}
\end{center}

{\em  So we obtain}
\begin{eqnarray*}
\P_3(10)&=&B(2,2;1,1)+B(0,0;3,3)+B(4,3;0,0)+B(2,1;2,2)+B(4,4;0,0)\\
&+&B(2,2;2,2)+B(4,2;1,1)+B(2,2;2,1)+B(0,0;4,3)+B(0,0;4,4)\\
&+&B(4,3;1,1)+B(6,3;0,0)+B(2,1;3,3)+B(2,1;3,2)+B(4,1;2,2)\\
&+&B(6,4;0,0)+B(4,4;1,1)+B(4,2;2,2)+B(2,2;3,2)+B(2,2;3,3)\\
&+&B(6,2;1,1)+B(4,2;2,1)+B(2,2;3,1)\\
&=&(1+1+1+1+1)+(2+2+1+1+1)+(2+3+1+1+1)\\
&+&(3+1+4+2+2)+(3+2+1)=38.
\end{eqnarray*}

\medskip
{\em 
On the other hand, by an exhaustive search, we find all equivalence classes of $t$-codes representing Type 3 alternating oriented pretzel links with $10$ crossings, listed below with one $t$-code from each equivalence class. 
Note that there are exactly 38 of them, as expected. With the understanding that these are all Type 3 $t$-codes, we have omitted $P_3$ from the codes for simplicity.

\medskip
\begin{tabular}{llll}
$ (2;2,2,2,2)$ &$ (3;2,2,3)$&$ (1;2,2,5)$&$ (1;2,3,4)$\\
$ (1;3,3,3)$&$ (0;2,2,2,4)$&$ (0;2,2,3,3)$&$ (0;2,3,2,3)$\\$ (4;2,2,-2)$&$ (2;2,2,-4)$&$ (2;2,4,-2)$&$ (2;3,3,-2)$\\
$ (1;2,2,3,-2)$&$ (1;2,3,2,-2)$&$ (0;2,2,2,2,-2)$&$ (0;2,2,-6)$\\
$ (0;2,4,-4)$&$ (0;2,6,-2)$&$ (0;3,3,-4)$&$ (0;4,4,-2)$\\
$ (0;3,5,-2)$&$ (3;3,-2,-2)$&$ (2;2,2,-2,-2)$&$ (2;2,-2,2,-2)$\\
$ (1;5,-2,-2)$&$ (0;2,4,-2,-2)$&$ (0;2,-2,4,-2)$&$ (0;3,3,-2,-2)$\\
$ (0;3,-2,3,-2)$&$ (1;3,-2,-4)$&$ (0;2,2,-2,-4)$&$ (0;2,-2,2,-4)$\\
$ (4;-2,-2,-2)$&$ (1;3,-2,-2,-2)$&$ (0;2,2,-2,-2,-2)$&$ (0;2,-2,2,-2,-2)$\\
$ (2;-2,-2,-4)$&$ (2;-2,-2,-2,-2)$&&
\end{tabular}
}
\end{example}

\section{Numerical Results}\label{sec_num}

\begin{figure}[h!]
\begin{center}
\includegraphics[scale=0.4]{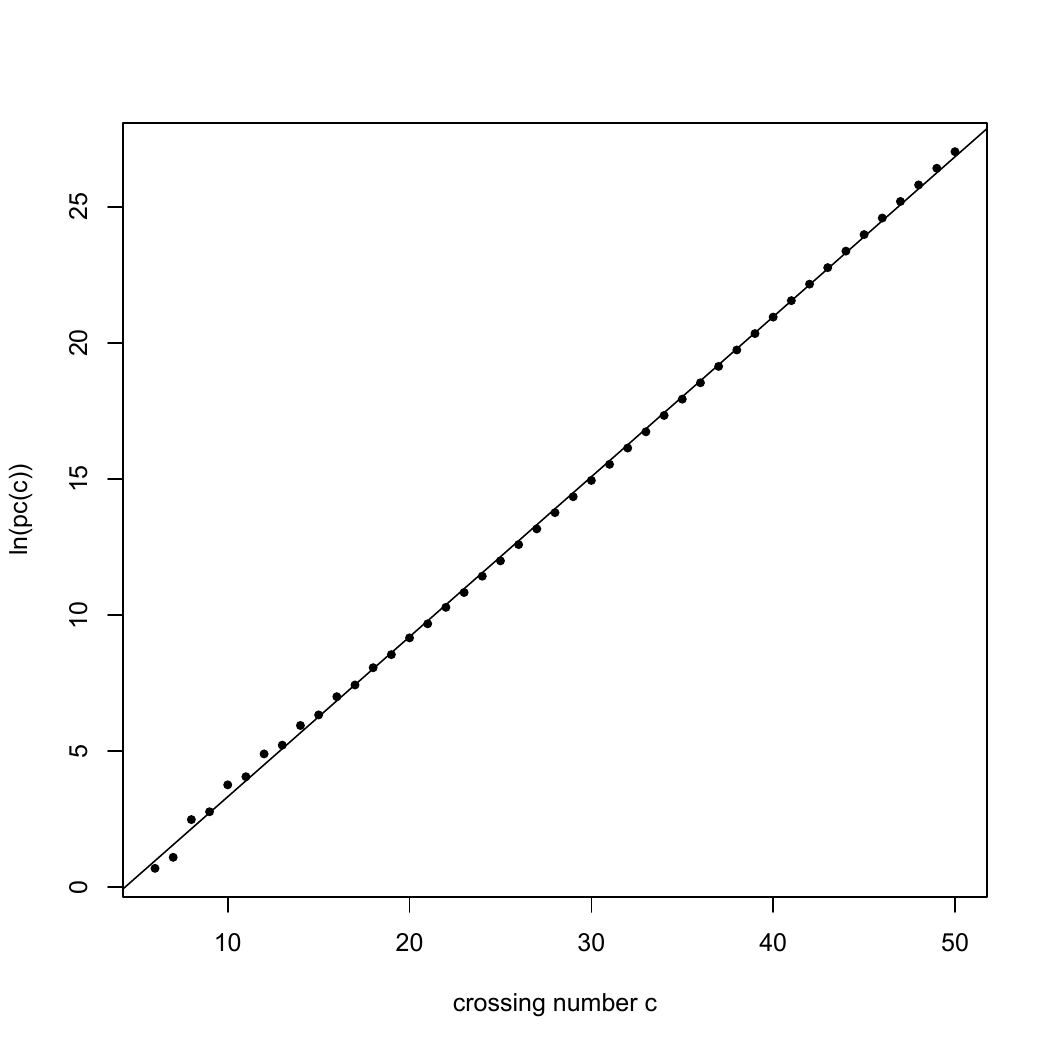}
\end{center}
\caption{The plot of $\ln(\P(c))$ shows a clear linear pattern. The fitting function for $\P(c)$ is $\P(c)\approx 0.0775e^{0.588c}$ with $r^2=0.9995$.
\label{data_plot}}
\end{figure}

Table \ref{data_table} below shows the numerical results for crossing numbers up to 50. Keep in mind that $\P_j(c)$ is the total number of mirror image pairs of Type $j$ alternating oriented pretzel links with crossing number $c$ ($1\le j\le 3$), and the total number of alternating oriented pretzel links with a given crossing number $c$ is $2\P(c)=2(\P_1(c)+\P_2(c)+\P_3(c))$. Table \ref{data_table} clearly shows that, as the crossing number increases, the number of Type 3 alternating oriented pretzel links overwhelmingly dominates the numbers of Type 1 and Type 2 alternating oriented pretzel links. 
Furthermore, $\P(c)$ demonstrates an exponential growth as a function of $c$ as shown in Figure \ref{data_plot}. 

\begin{center}
\begin{longtable}{|c||c|c|c||c|}
\caption{The numbers of Type 1, 2 and 3 alternating oriented pretzel link mirror image pairs  with a given crossing number are listed in columns $\P_1(c)$,  $\P_2(c)$ and $\P_3(c)$. The last column $\P(c)$ lists the summation of $\P_1(c)$,  $\P_2(c)$ and $\P_3(c)$.}
\label{data_table}\\
\hline $ c $ & $\P_1(c)$ & $ \P_2(c) $   & $ \P_3(c) $ & $ \P(c)$\\
\hline\hline
$ \le 5 $ & $ 0$ & $ 0 $   & $ 0$ & $ 0$ 
\\
\hline
$ 6 $ & $ 0$ & $ 1 $   & $ 1$ & $ 2$ 
\\
\hline
$ 7 $ & $ 0$ & $ 0 $   & $ 3$ & $ 3$ 
\\
\hline
$ 8 $ & $ 0$ & $ 2 $   & $ 10$ & $ 12$ 
\\
\hline
$ 9 $ & $ 1$ & $ 0 $   & $ 15$ & $ 16$ 
\\
\hline
$ 10 $ & $ 1$ & $ 4 $   & $ 38$ & $ 43$ 
\\
\hline
$ 11 $ & $ 2$ & $ 0 $   & $ 56$ & $ 58$ 
\\
\hline
$ 12 $ & $ 3$ & $ 8 $   & $ 123$ & $ 134$ 
\\
\hline
$ 13 $ & $ 5$ & $ 0 $   & $ 180$ & $ 185$ 
\\
\hline
$ 14 $ & $ 6$ & $ 13 $   & $ 362$ & $ 381$ 
\\
\hline
$ 15 $ & $ 11$ & $ 0 $   & $ 551$ & $ 562$ 
\\
\hline
$ 16 $ & $ 14$ & $ 24 $   & $ 1060$ & $ 1098$ 
\\
\hline
$ 17 $ & $ 20$ & $ 0 $   & $ 1670$ & $ 1690$ 
\\
\hline
$ 18 $ & $ 26$ & $ 40 $   & $ 3122$ & $ 3188$ 
\\
\hline
$ 19 $ & $ 36$ & $ 0 $   & $ 5122$ & $ 5158$ 
\\
\hline
$ 20 $ & $ 47$ & $ 71 $   & $ 9426$ & $ 9544$ 
\\
\hline
$ 21 $ & $ 65$ & $ 0 $   & $ 15947$ & $ 16012$ 
\\
\hline
$ 22 $ & $ 83$ & $ 119 $   & $ 29099$ & $ 29301$ 
\\
\hline
$ 23 $ & $ 110$ & $ 0 $   & $ 50429$ & $ 50539$ 
\\
\hline
$ 24 $ & $ 143$ & $ 216 $   & $ 91701$ & $ 92060$ 
\\
\hline
$ 25 $ & $ 188$ & $ 0 $   & $ 161588$ & $ 161776$ 
\\
\hline
$ 26 $ & $ 241$ & $ 372 $   & $293479 $ & $ 294092$ 
\\
\hline
$ 27 $ & $ 315$ & $ 0 $   & $ 523293$ & $ 523608$ 
\\
\hline
$ 28 $ & $ 405$ & $ 678 $   & $ 950725$ & $ 951808$ 
\\
\hline
$ 29 $ & $ 524$ & $ 0 $   & $ 1708860$ & $ 1709384$ 
\\
\hline
$ 30 $ & $675$ & $ 1215 $   & $ 3107773$ & $3109663$ 
\\
\hline
$ 31 $ & $ 871$ & $ 0 $   & $ 5617881$ & $5618752$ 
\\
\hline
$ 32 $ & $ 1120$ & $ 2240 $   & $ 10230499$ & $10233859$ 
\\
\hline
$ 33 $ & $ 1446$ & $ 0 $   & $18569086$ & $18570532$ 
\\
\hline
$ 34 $ & $ 1859$ & $ 4102 $   & $ 33864326$ & $33870287$ 
\\
\hline
$ 35 $ & $ 2379$ & $ 0 $   & $ 61654785$ & $61657182$ 
\\
\hline
$ 36 $ & $ 3088$ & $ 7674 $   & $ 112602737$ & $112613499$ 
\\
\hline
$ 37 $ & $ 3979$ & $ 0 $   & $ 205497471$ & $ 205501450$ 
\\
\hline
$ 38 $ & $ 5126$ & $ 14299 $   & $ 375831251$ & $ 375850676$ 
\\
\hline
$ 39 $ & $ 6613$ & $ 0 $   & $ 687206188$ & $687212801$ 
\\
\hline
$ 40 $ & $ 8531$ & $ 27000 $   & $ 1258468810$ & $1258504341$ 
\\
\hline
$ 41 $ & $ 11009$ & $ 0 $   & $ 2304807470$ & $ 2304818479$ 
\\
\hline
$ 42 $ & $ 14217$ & $ 50952 $   & $4225898392$ & $ 4225963561$ 
\\
\hline
$ 43 $ & $ 18364$ & $ 0 $   & $7750154298$ & $ 7750172662$ 
\\
\hline
$ 44 $ & $ 23736$ & $ 96896 $   & $ 14226040972$ & $  14226161604$ 
\\
\hline
$ 45 $ & $ 30696$ & $ 0 $   & $ 26121609372$ & $ 26121640068$ 
\\
\hline
$ 46 $ & $ 39713$ & $ 184397 $   & $  47998211056$ & $ 47998435166$ 
\\
\hline
$ 47 $ & $ 51399$ & $ 0 $   & $  88228272471$ & $  88228323870$ 
\\
\hline
$ 48 $ & $ 66571$ & $ 352684 $   & $ 162274113329$ & $ 162274532584$ 
\\
\hline
$ 49 $ & $ 86243$ & $ 0 $   & $ 298574262099$ & $ 298574348342$ 
\\
\hline
$ 50 $ & $ 111794$ & $ 675174 $   & $ 549639730670$ & $ 549640517638$ 
\\
\hline
\end{longtable}
\end{center}

Let us compare $\P(c)$ with the number of oriented rational links with crossing number $c$, which grows at a rate of $\frac{2^c}{6}\approx .167e^{.693c}$. 
We see that $\P(c)$ grows at a slower rate when compared to the number of oriented rational links with crossing number $c$. 
However in this paper we only considered the alternating oriented pretzel links. 
It is intuitive that for large crossing numbers, there would be many more non-alternating oriented pretzel links than the alternating ones. The enumeration of the non alternating oriented pretzel links would be a more challenging task which we shall tackle in the future. 
So it is possible that the number of all pretzel links with a given crossing number may surpass the number of rational links for large crossing numbers. The enumeration results obtained here show us that the set of pretzel links is large. 

\medskip
We shall end our paper with the following remark.

\medskip
\begin{remark}{\em
The precise formula for the braid index of any rational link has been known for a long time \cite{Murasugi1991}. 
The derivation of this formula is based on a specific non-minimal diagram of the rational link.
In a more recent work, it has been shown that a formulation of the braid index can in fact be read out directly from any minimum diagram of the rational link \cite{Diao2021_1}. 
In the same paper \cite{Diao2021_1}, it has been shown that a formula for the braid index of any alternating oriented pretzel link can also be derived from any minimum diagram of the link. 
However, the situation for non alternating oriented pretzel links is quite different. 
Some non alternating oriented pretzel links possess different properties and the same techniques that work for the alternating oriented pretzel links and rational links (which are all alternating) fail on them \cite{Diao2024_1,Diao2024_2}. 
This makes the set of non alternating oriented pretzel links a testing bed for future tools aimed at computing the braid index.
}
\end{remark}

\section*{Acknowledgement}
Charlotte Aspinwall and Tobias Clark are grateful for the support from the 2024 UNC Charlotte summer undergraduate research program funded by NSF-REU DMS-2150179.

\section{Appendix}

The proof of Formula (\ref{hformula1}), that is,  if $k=k_1+k_2$ is odd then
$$
h(n_1,k_1;n_2,k_2)
=
\begin{cases}
 \frac{1}{2}\binom{\frac{k-1}{2}}{\frac{k_2}{2}}\binom{\frac{n_1-1}{2}}{\frac{k_1-1}{2}}\binom{\frac{n_2}{2}-1}{\frac{k_2}{2}-1}, & \ {\rm if}\ n_1\equiv k_1\equiv 1\ {\rm mod}\ 2,\ n_2\equiv k_2\equiv 0 \ {\rm mod}\ 2,\\
 \frac{1}{2}\binom{\frac{k-1}{2}}{\frac{k_1}{2}}\binom{\frac{n_1}{2}-1}{\frac{k_1}{2}-1}\binom{\frac{n_2-1}{2}}{\frac{k_2-1}{2}}, & \ {\rm if}\ n_1\equiv k_1\equiv 0\ {\rm mod}\ 2,\ n_2\equiv k_2\equiv 1 \ {\rm mod}\ 2,\\
 \frac{1}{2}\binom{\frac{k-1}{2}}{\frac{k_1}{2}}\binom{\frac{n_1}{2}-1}{\frac{k_1}{2}-1}\binom{\frac{n_2}{2}-1}{\frac{k_2-1}{2}}, & \ {\rm if}\ n_1\equiv k_1\equiv n_2\equiv 0\ {\rm mod}\ 2,\ k_2\equiv 1 \ {\rm mod}\ 2,\\
  \frac{1}{2}\binom{\frac{k-1}{2}}{\frac{k_2}{2}}\binom{\frac{n_1}{2}-1}{\frac{k_1-1}{2}}\binom{\frac{n_2}{2}-1}{\frac{k_2}{2}-1}, & \ {\rm if}\ n_1\equiv k_2\equiv n_2\equiv 0\ {\rm mod}\ 2,\ k_1\equiv 1 \ {\rm mod}\ 2,\\
 0, & \ {\rm if}\  k_1\equiv  n_2\equiv  1\ {\rm mod}\ 2\ {\rm or}\  k_2\equiv  n_1\equiv  1\ {\rm mod}\ 2.
  \end{cases}
$$

Let us keep in mind that by (\ref{B_{k}(t,s)}), $h(n_1,k_1;n_2,k_2)$ is the coefficient of the $t^{n_1}s^{n_2}$ power term contributed by $k_1$ copies of $f(t)$ and $k_2$ copies of $f(s)$ from
$$
\frac{1}{2} \left(f(t)+f(s)\right)\left(f(t^2)+f(s^2)\right)^{(k-1)/2}=\frac{1}{2} \left(\frac{t}{1-t}+\frac{s}{1-s}\right)\left(\frac{t^2}{1-t^2}+\frac{s^2}{1-s^2}\right)^{(k-1)/2}.
$$
Thus if $k_1$ is odd, then only $\frac{1}{2}f(t)\left(f(t^2)+f(s^2)\right)^{(k-1)/2}$ can contribute $k_1$ copies of $f(t)$ and $k_2$ copies of $f(s)$. However in the power expansion of $\frac{1}{2}f(t)\left(f(t^2)+f(s^2)\right)^{(k-1)/2}$ the power of $s$ in every term is even, therefore there will be no power term of the form $t^{n_1}s^{n_2}$ if $n_2$ is odd. That is, $h(n_1,k_1;n_2,k_2)=0$ if $k_1$ and $n_2$ are both odd. Same can be said if $k_2$ and $n_1$ are both odd. This proves the last case in (\ref{hformula1}).

\medskip
Let us now consider the case $n_1\equiv k_1\equiv 1\ {\rm mod}\ 2,\ n_2\equiv k_2\equiv 0 \ {\rm mod}\ 2$. As we mentioned above, only $\frac{1}{2}f(t)\left(f(t^2)+f(s^2)\right)^{(k-1)/2}$ can contribute $k_1$ copies of $f(t)$ and $k_2$ copies of $f(s)$ if $k_1$ is odd, and the combination of these copies is
\begin{equation}\label{append_eq1}
\frac{1}{2}\binom{\frac{k-1}{2}}{\frac{k_2}{2}}f(t)\left(f(t^2)\right)^{\frac{k_1-1}{2}}\left(f(s^2)\right)^{\frac{k_2}{2}}=\frac{1}{2}\binom{\frac{k-1}{2}}{\frac{k_2}{2}}\left(\frac{t}{1-t}\right)\left(\frac{t^2}{1-t^2}\right)^{\frac{k_1-1}{2}}\left(\frac{s^2}{1-s^2}\right)^{\frac{k_2}{2}}.
\end{equation}
We then extract the $s^{n_2}$ power term from $\left(\frac{s^2}{1-s^2}\right)^{\frac{k_2}{2}}$, which is $\binom{n_2/2-1}{k_2/2-1}s^{n_2}$ by (\ref{MacSeries}), and extract the $t^{n_1}$ power term from  
\begin{eqnarray}
&&\left(\frac{t}{1-t}\right)\left(\frac{t^2}{1-t^2}\right)^{\frac{k_1-1}{2}}\label{append_eq2}\\
&=&t^{k_1}\left(1+t+t^2+\cdots \right)\left(1+\frac{k_1-1}{2}t^2+\cdots+\binom{\frac{k_1-1}{2}+j-1}{j}t^{2j}+\cdots \right)\nonumber,
\end{eqnarray}
which is 
$$
\left(\sum_{j=0}^{\frac{n_1-k_1}{2}}\binom{\frac{k_1-1}{2}+j-1}{j}\right)t^{n_1}=\left(\sum_{j=0}^{\frac{n_1-k_1}{2}}\binom{\frac{k_1-1}{2}+j-1}{\frac{k_1-1}{2}-1}\right)t^{n_1}.
$$
By (\ref{binom_iden_1}), we have
$$
\sum_{j=0}^{\frac{n_1-k_1}{2}}\binom{\frac{k_1-1}{2}+j-1}{\frac{k_1-1}{2}-1}=\binom{\frac{n_1-1}{2}}{\frac{k_1-1}{2}}.
$$
Thus, the coefficient of the $t^{n_1}s^{n_2}$ power term is
$$
\frac{1}{2}\binom{\frac{k-1}{2}}{\frac{k_2}{2}}\binom{\frac{n_1-1}{2}}{\frac{k_1-1}{2}}\binom{\frac{n_2}{2}-1}{\frac{k_2}{2}-1}.
$$
This completes the proof for this case. The case $n_1\equiv k_1\equiv 0\ {\rm mod}\ 2,\ n_2\equiv k_2\equiv 1 \ {\rm mod}\ 2$ can be similarly proved by simply switching the roles of $t$ and $s$.

\smallskip
Finally let us consider the case $n_1\equiv k_2\equiv n_2\equiv 0\ {\rm mod}\ 2,\ k_1\equiv 1 \ {\rm mod}\ 2$. Similar to the previous case, we are to extract the $t^{n_1}s^{n_2}$ power term from (\ref{append_eq1}). Here we will again extract the $s^{n_2}$ power term from $\left(\frac{s^2}{1-s^2}\right)^{\frac{k_2}{2}}$, which is $\binom{n_2/2-1}{k_2/2-1}s^{n_2}$, and extract the $t^{n_1}$ power term  from (\ref{append_eq2}), whose coefficient is
$$
\left(\sum_{j=0}^{\frac{n_1-k_1-1}{2}}\binom{\frac{k_1-1}{2}+j-1}{\frac{k_1-1}{2}-1}\right)=\binom{\frac{n_1}{2}-1}{\frac{k_1-1}{2}}
$$
by (\ref{binom_iden_1}). Thus the coefficient of $t^{n_1}s^{n_2}$ is
$$
\frac{1}{2}\binom{\frac{k-1}{2}}{\frac{k_2}{2}}\binom{\frac{n_1}{2}-1}{\frac{k_1-1}{2}}\binom{\frac{n_2}{2}-1}{\frac{k_2}{2}-1}.
$$ 
This completes the proof of Formula (\ref{hformula1}) since the remaining case $n_1\equiv k_1\equiv n_2\equiv 0\ {\rm mod}\ 2,\ k_2\equiv 1 \ {\rm mod}\ 2$ is symmetric to what we just proved.

\medskip
The proof of Formula (\ref{hformula2}), that is,  if $k$ is even but $k_1$, $k_2$ are both odd, then
$$
h(n_1,k_1;n_2,k_2)
=
\begin{cases}
 \frac{1}{2}\binom{\frac{k}{2}-1}{\frac{k_1-1}{2}}\binom{\frac{n_1-1}{2}}{\frac{k_1-1}{2}}\binom{\frac{n_2-1}{2}}{\frac{k_2-1}{2}}, & \ {\rm if}\ n_1\equiv n_2\equiv 1\ {\rm mod}\ 2,\\
 \frac{1}{2}\binom{\frac{k}{2}-1}{\frac{k_1-1}{2}}\binom{\frac{n_1}{2}-1}{\frac{k_1-1}{2}}\binom{\frac{n_2-1}{2}}{\frac{k_2-1}{2}}, & \ {\rm if}\ n_1\equiv n_2+1\equiv 0\ {\rm mod}\ 2,\\
 \frac{1}{2}\binom{\frac{k}{2}-1}{\frac{k_1-1}{2}}\binom{\frac{n_1-1}{2}}{\frac{k_1-1}{2}}\binom{\frac{n_2}{2}-1}{\frac{k_2-1}{2}}, & \ {\rm if}\ n_1\equiv n_2+1\equiv 1\ {\rm mod}\ 2,\\
 \frac{1}{2}\binom{\frac{k}{2}-1}{\frac{k_1-1}{2}}\binom{\frac{n_1}{2}-1}{\frac{k_1-1}{2}}\binom{\frac{n_2}{2}-1}{\frac{k_2-1}{2}}, & \ {\rm if}\ n_1\equiv n_2\equiv 0\ {\rm mod}\ 2.
\end{cases}
$$

\medskip
By (\ref{B_{k}(t,s)}), $h(n_1,k_1;n_2,k_2)$ is the coefficient of the $t^{n_1}s^{n_2}$ power term contributed by $k_1$ copies of $f(t)$ and $k_2$ copies of $f(s)$ from
\begin{eqnarray*}
&&\frac{1}{4} \left(f(t)+f(s)\right)^2\left(f(t^2)+f(s^2)\right)^{\frac{k-2}{2}}+\frac{1}{4}\left(f(t^2)+f(s^2)\right)^{\frac{k}{2}}\\
&=&
\frac{1}{4}\left(f^2(t)+f^2(s)\right)\left(f(t^2)+f(s^2)\right)^{\frac{k-2}{2}}+\frac{1}{2} f(t)f(s)\left(f(t^2)+f(s^2)\right)^{\frac{k-2}{2}}+\frac{1}{4}\left(f(t^2)+f(s^2)\right)^{\frac{k}{2}}.
\end{eqnarray*}
Since $\frac{1}{4}\left(f(t^2)+f(s^2)\right)^{\frac{k}{2}}$, $\frac{1}{4}f^2(t)\left(f(t^2)+f(s^2)\right)^{\frac{k-2}{2}}$ and $\frac{1}{4}f^2(s)\left(f(t^2)+f(s^2)\right)^{\frac{k-2}{2}}$ can only contribute even copies of $f(t)$ and $f(s)$, we only need to consider the term $\frac{1}{2} f(t)f(s)\left(f(t^2)+f(s^2)\right)^{\frac{k-2}{2}}$. The combined terms with $k_1$ copies of $f(t)$ and $k_2$ copies of $f(s)$ from it is
\begin{eqnarray*}
&&\frac{1}{2}\binom{\frac{k}{2}-1}{\frac{k_1-1}{2}}f(t)f(s)f^{\frac{k_1-1}{2}}(t^2)f^{\frac{k_2-1}{2}}(s^2)\\
&=&
\frac{1}{2}\binom{\frac{k}{2}-1}{\frac{k_1-1}{2}}\left(\frac{t}{1-t}\right)\left(\frac{t^2}{1-t^2}\right)^{\frac{k_1-1}{2}}\left(\frac{s}{1-s}\right)\left(\frac{s^2}{1-s^2}\right)^{\frac{k_2-1}{2}}.
\end{eqnarray*}
Now we just need to extract the $t^{n_1}$ term from 
$$
f(t)f^{\frac{k_1-1}{2}}(t^2)=\left(\frac{t}{1-t}\right)\left(\frac{t^2}{1-t^2}\right)^{\frac{k_1-1}{2}}
$$ 
and the $s^{n_2}$ term from 
$$
f(s)f^{\frac{k_2-1}{2}}(s^2)=\left(\frac{s}{1-s}\right)\left(\frac{s^2}{1-s^2}\right)^{\frac{k_2-1}{2}},
$$
 much like what we did in the proof of Formula (\ref{hformula1}) using the expansion formula given in (\ref{append_eq2}). 

\medskip
The proof of Formula (\ref{hformula3}), that is,  if $k_1$, $k_2$ are both even, then 
$$
h(n_1,k_1;n_2,k_2)
=
\begin{cases}
 \phantom{+}\frac{1}{4}\binom{\frac{k}{2}-1}{\frac{k_2}{2}}\binom{\frac{n_2}{2}-1}{\frac{k_2}{2}-1}\left(\binom{\frac{n_1}{2}}{\frac{k_1}{2}}+\binom{\frac{n_1}{2}-1}{\frac{k_1}{2}}\right) & \\
 +  \frac{1}{4}\binom{\frac{k}{2}-1}{\frac{k_1}{2}}\binom{\frac{n_1}{2}-1}{\frac{k_1}{2}-1}
 \left(\binom{\frac{n_2}{2}}{\frac{k_2}{2}}+\binom{\frac{n_2}{2}-1}{\frac{k_2}{2}}\right) &  \\
 + \frac{1}{4}\binom{\frac{k}{2}}{\frac{k_1}{2}}\binom{\frac{n_1}{2}-1}{\frac{k_1}{2}-1}\binom{\frac{n_2}{2}-1}{\frac{k_2}{2}-1}, & \ {\rm if}\ n_1\equiv n_2\equiv 0\ {\rm mod}\ 2,\\
\frac{1}{2}\binom{\frac{k}{2}-1}{\frac{k_1}{2}}\binom{\frac{n_1}{2}-1}{\frac{k_1}{2}-1}\binom{\frac{n_2-1}{2}}{\frac{k_2}{2}}, & \ {\rm if}\ n_1\equiv n_2+1\equiv 0\ {\rm mod}\ 2,\\ 
\frac{1}{2}\binom{\frac{k}{2}-1}{\frac{k_2}{2}}\binom{\frac{n_2}{2}-1}{\frac{k_2}{2}-1}\binom{\frac{n_1-1}{2}}{\frac{k_1}{2}}, & \ {\rm if}\ n_1+1\equiv n_2\equiv 0\ {\rm mod}\ 2,\\
0, &\ {\rm if}\ n_1\equiv n_2\equiv 1\ {\rm mod}\ 2.\\
\end{cases}
$$

By (\ref{B_{k}(t,s)}), $h(n_1,k_1;n_2,k_2)$ is the coefficient of the $t^{n_1}s^{n_2}$ power term contributed by $k_1$ copies of $f(t)$ and $k_2$ copies of $f(s)$ from
\begin{eqnarray*}
&&\frac{1}{4} \left(f(t)+f(s)\right)^2\left(f(t^2)+f(s^2)\right)^{\frac{k-2}{2}}+\frac{1}{4}\left(f(t^2)+f(s^2)\right)^{\frac{k}{2}}\\
&=&
\frac{1}{4}\left(f^2(t)+f^2(s)\right)\left(f(t^2)+f(s^2)\right)^{\frac{k-2}{2}}+\frac{1}{2} f(t)f(s)\left(f(t^2)+f(s^2)\right)^{\frac{k-2}{2}}+\frac{1}{4}\left(f(t^2)+f(s^2)\right)^{\frac{k}{2}}.
\end{eqnarray*}
We can disregard the term $\frac{1}{2} f(t)f(s)\left(f(t^2)+f(s^2)\right)^{\frac{k-2}{2}}$ in the above since it can only contribute terms containing odd copies of $f(t)$ and $f(s)$. Since the expansions of $f(t^2)=t^2/(1-t^2)$ and $f(s^2)=s^2/(1-s^2)$ contain only even powers of $t$ and $s$ respectively, out of the three terms $\frac{1}{4} f^2(t)\left(f(t^2)+f(s^2)\right)^{\frac{k-2}{2}}$, $\frac{1}{4} f^2(s)\left(f(t^2)+f(s^2)\right)^{\frac{k-2}{2}}$ and $\frac{1}{4}\left(f(t^2)+f(s^2)\right)^{\frac{k}{2}}$, only $\frac{1}{4} f^2(t)\left(f(t^2)+f(s^2)\right)^{\frac{k-2}{2}}$ contains odd powers of $t$ and only $\frac{1}{4} f^2(s)\left(f(t^2)+f(s^2)\right)^{\frac{k-2}{2}}$ contains odd power of $s$. It follows that $h(n_1,k_1;n_2,k_2)=0$ if $n_1$ and $n_2$ are both odd. The case that $n_1$ and $n_2$ are both even have been proved in Section \ref{sec_enum}. So the only remaining cases to be proven are the two symmetric cases when exactly one of $n_1$, $n_2$ is odd. Let us prove the case $n_1$ is odd and $n_2$ is even. In this case we just need to extract the $t^{n_1}s^{n_2}$ term from 
\begin{equation}\label{append_eq2}
\frac{1}{4}\binom{\frac{k}{2}-1}{\frac{k_2}{2}}\left(\frac{t}{1-t}\right)^2\left(\frac{t^2}{1-t^2}\right)^{\frac{k_1}{2}-1}\left(\frac{s^2}{1-s^2}\right)^{\frac{k_s}{2}},
\end{equation}
which is the combined terms in $\frac{1}{4} f^2(t)\left(f(t^2)+f(s^2)\right)^{\frac{k-2}{2}}$ with $k_1$ copies of $f(t)$ and $k_2$ copies of $f(s)$. As we have seen in Section \ref{sec_enum}, the $s^{n_2}$ term contributed from $\left(\frac{s^2}{1-s^2}\right)^{\frac{k_2}{2}}$ is of the form 
\begin{equation}\label{append_eq5}
\binom{\frac{n_2}{2}-1}{\frac{k_2}{2}-1}s^{n_2}.
\end{equation}
 By (\ref{expans_k1even}), we have
\begin{eqnarray}
&&\left(\frac{t}{1-t}\right)^2\left(\frac{t^2}{1-t^2}\right)^{\frac{k_1}{2}-1}\label{append_eq3}\\
&=&t^{k_1}\left(1+2t+3t^2+\cdots\right)\left(1+\binom{\frac{k_1}{2}-2+1}{\frac{k_1}{2}-2}t^2+\cdots+\binom{\frac{k_1}{2}-2+j}{\frac{k_1}{2}-2}t^{2j}+\cdots\right).\nonumber
\end{eqnarray}
So the coefficient of the $t^{n_1}$ term extracted from it is
\begin{equation}\label{append_eq4}
\sum_{j=1}^{\frac{n_1-k_1+1}{2}}2j\binom{\frac{n_1-1}{2}-j-1}{\frac{k_1}{2}-2}=2\sum_{j=1}^{\frac{n_1-k_1+1}{2}}j\binom{\frac{n_1-1}{2}-j-1}{\frac{k_1}{2}-2}=2\binom{\frac{n_1-1}{2}}{\frac{k_1}{2}}
\end{equation}
by (\ref{binom_iden_3}). Notice that the above derivation is only valid under the assumption $k_1\ge 4$. In the case that $k_1=2$, the $t^{n_1}$ term is extracted from
$t^{2}(1+2t+3t^2+\cdots)$, which has coefficient $n_1-1=2\binom{\frac{n_1-1}{2}}{1}$, so (\ref{append_eq4}) holds for $k_1=2$ as well. Combining (\ref{append_eq2}), (\ref{append_eq5}) and (\ref{append_eq4}), we obtain
$$
h(n_1,k_1;n_2,k_2)=\frac{1}{2}\binom{\frac{k}{2}-1}{\frac{k_2}{2}}\binom{\frac{n_2}{2}-1}{\frac{k_2}{2}-1}\binom{\frac{n_1-1}{2}}{\frac{k_1}{2}},
$$
as desired.
\end{document}